%% file: jcpE1.tex
\documentclass[reviewcopy]{elsart}

\usepackage{amssymb}
\usepackage{amsmath}
\usepackage[usenames,dvips]{color}
\usepackage{lineno}
\usepackage{graphicx}
\usepackage{bm}

\def\be{\begin{equation}}
\def\ee{\end{equation}}

\def\ba{\begin{array}{lllll}}
\def\ea{\end{array}}

\def\dsp{\displaystyle}
 
\def\mathbi#1{{\boldsymbol{#1}}}

\def\bu{\mathbi{u}}
\def\bv{\mathbi{v}}
\def\bn{\mathbi{n}}
\def\bw{\mathbi{w}}
\def\bff{\mathbi{f}}

\def\bx{\mathbi{x}}

\def\buref{\mathbi{u}_{{\rm ref}}}
\def\Tref{{T}_{{\rm ref}}}
\def\psiref{{\psi}_{{\rm ref}}}
\def\pref{{p}_{{\rm ref}}}

\def\points{{\cal P}}
\def\clusters{{\cal G}}

\def\cv{K}
\def\cvv{L}
\def\cvvv{M}
\def\cvvvv{N}

\def\d{{\rm d}}
\def\dcvedge{d_{\cv,\edge}}

\def\dcvvedge{d_{\cvv,\edge}}

\def\disc{{\mathcal D}}

\def\dcvedge{d_{\cv,\edge}}

\def\dr{\partial}

\def\edge{{\sigma}}
\def\edgep{{\sigma'}}
\def\edges{{\cal E}}
\def\edgesint{{\cal E}_{\mathrm{int}}}
\def\edgesext{{\cal E}_{\mathrm{ext}}}

\def\edgescv{{\cal E}_\cv}

\def\edgescvv{{\cal E}_\cvv}

\def\frontomega{{\Gamma} }
\def\frontomegaun{{\Gamma_1}}
\def\frontomegadeux{{\Gamma_2}}

\def\hmesh{{\cal H}_{\mesh}(\O)}

\def\lap{\Delta}

\def\mcv{\meas_\cv}
\def\mcvv{\meas_\cvv}
\def\meas{{\rm m}}
\def\medge{{\rm m}_{\edge}}
\def\medgep{{\rm m}_{\edge'}}

\def\mesh{{\cal M}}

\def\mK{{\mathrm m}_K}

\def\ncvedge{\bn_{\cv,\edge}}
\def\ncvedgep{\bn_{\cv,\edge'}}
\def\ncvedgep{\bn_{\cv,\edge'}}
\def\ncvvedge{\bn_{\cvv,\edge}}

\def\normedeu(#1){\|#1\|_{2}}
\def\normeinf(#1){\|#1\|_{\infty}}
\def\normeH1(#1){\|#1\|_{H_1}}

\def\O{\Omega}

\def\R{\mathbb{R}}

\def\xcv{{\bx}_\cv}
\def\xcvv{{\bx}_\cvv}

\def\xedge{{\bx}_\edge}
\def\xedgep{{\bx}_\edge'}

\def\xedgep{{\bx}_{\edge'}}

\def\b0{{\mathbi{0}}}
\def\Ra{{\rm{R\!a}}}
\def\Pr{{\rm{P\!r}}}

\def\Nu{{\rm{N\!u}}}

\def\ex{{\mathbi{e}_1}}

\def\ez{{\mathbi{e}_3}}
\def\ei{{\mathbi{e}_i}}
\def\bn{{\mathbi{n}}}
\def\eex{{\mathbi{e}_x}}
\def\eey{{\mathbi{e}_y}}
\def\0b{{O\!\!\!\!/}}

\def\dset{X^\disc}

\def\dsetU{X_{0}^\disc}
\def\dsetTheta{X_{\Gamma_1,0}^\disc}

\def\xs{{\bx}_{{\rm s}}}
\def\xsi{{\bx}_{{\rm s}}}

\def\xsi(#1){\bx^{(#1)}_{{\rm s}}}
\def\xcvi(#1){\bx_\cv^{(#1)}}
\def\ui(#1){u^{(#1)}}
\def\vi(#1){v^{(#1)}}

\def\div{{\rm div}}
\def\bdiv{{\mbox {\bf div}}}

\def\grad{\boldsymbol{\nabla}}

\def\tstf{{v}}
 
\begin{document}

\begin{frontmatter}





\title{A collocated finite volume scheme to solve free convection for general non-conforming grids}


\author{Eric Ch\'enier\corauthref{cor1}},
\ead{Eric.Chenier@univ-paris-est.fr}
\corauth[cor1]{Corresponding author.}
\author{Robert Eymard}
\ead{Robert.Eymard@univ-paris-est.fr}
\address{Universit\'e Paris-Est, Laboratoire Mod\'elisation et Simulation Multi Echelle, 
MSME FRE3160 CNRS, 5 bd Descartes, 77454 Marne-la-Vall\'ee, France}
\author{Rapha\`ele Herbin}
\ead{Raphaele.Herbin@latp.univ-mrs.fr}
\address{Universit\'e de Provence, 39 rue Joliot-Curie, 13453 Marseille, France}

\begin{abstract}
We present a new collocated numerical scheme for the approximation of the Navier-Stokes and energy equations under the Boussinesq assumption for general grids, using the velocity-pressure unknowns. 
This scheme is based on a recent scheme for the diffusion terms.
Stability properties are drawn from particular choices for the pressure gradient and the non-linear terms. 
Convergence of the approximate solutions may be proven mathematically.
Numerical results show the accuracy of the scheme on irregular grids.
\end{abstract}

\begin{keyword}
Collocated finite volume schemes \sep general non-conforming grids \sep Navier-Stokes equations \sep Boussinesq assumption

\PACS 65C20 \sep 76R10 
\end{keyword}
\end{frontmatter}

\section{Introduction}
\label{Introduction}

Finite volume methods have been widely used in computational fluid dynamics for a long time; they are well adapted to the discretization of partial differential equations under conservative form;
one of their attractive features is that the resulting discretized equation has a clear physical interpretation  \cite{Patankar1980}. 
In the framework of incompressible fluid flows, two strategies are often opposed, namely staggered and collocated schemes.
The staggered strategy, which has become very popular since Patankar's book \cite{Patankar1980},   
remains mainly restricted to geometrical domains with parallel and orthogonal boundary faces. 
Therefore, for computations on complex domains with general meshes, the collocated strategy which consists in approximating all unknowns on the same set of points (called collocation points but also cell centers or simply centers), is often preferred, even though the pressure-velocity coupling demands some cure for the stabilization of the well-known checkerboard pressure modes; to this purpose, various pressure stabilization procedures, based on improvements of the {\em Momentum Interpolation Method}  proposed by Rhie and Chow \cite{Rhie1983}, are frequently used \cite{Nagele2007}.\\
In \cite{Chenier05,Touazi2007}, a collocated finite volume scheme for incompressible flows is developed on so called ``admissible" unstructured meshes, that are meshes satisfying the two following conditions: the straight line joining the centers of two adjacent control volumes is perpendicular to the common edge, and the neighboring control volumes and the associated centers are arranged in the same order, with respect to the common edge. 
Rectangular or orthogonal parallelepipedic meshes, triangular (2D) or tetrahedral (3D) Delaunay meshes, and Voronoi meshes fulfill these requirements.
Under this assumption, the isotropic diffusion fluxes can be consistently approximated by a two-point finite difference scheme.
Using this approximation for a pure diffusion problem yields a symmetric ``M-matrix'' (which ensures monotony); the  stencil is limited to the control volume itself and its natural neighbors and it leads to the classical 5- and 7-point schemes on rectangles and orthogonal parallelepipeds.
Unfortunately, although the use of such grids considerably widens the variety of geometric shapes which can be gridded, it is far from solving all the critical needs resulting from actual problems:
\begin{itemize}
\item For complex 3D domains, it is well known that the use of a  large number of flat tetrahedra produces high discretization errors; generalized hexahedric meshes are often preferred: these are made of 3D elevations of quadrangular meshes, for which the faces of the control volumes may no longer be planar.
\item To our knowledge, there is yet no mesh software able to grid any geometrical shape in 3D
using Voronoi or Delaunay tessellations while respecting the boundaries and the local refinement requirements.
\item In compressible flows, the approximation of the full tensor by the usual two point scheme is no longer consistent even on admissible meshes; multi-point approximations are therefore required. 
\item Boundary layers are classically meshed with refined grids, so that the discretization scheme should be able to deal with non-conforming meshes.
\end{itemize}
Whereas there is no real difficulty to discretize the convective terms for general  non-conforming grids, writing accurate diffusion approximations, particularly relevant for low Reynolds (P\'eclet) flows, is still a challenge on such meshes.\\
In the early 80's, Kershaw \cite{Kershaw1981} first proposed a nine-point scheme on structured quadrilateral grids by using the restrictive assumption of a smooth mapping between the logical mesh and the spatial coordinates.
Since then, numerous works have been published to efficiently solve the diffusion equations in general geometry (see \cite{Breil2007} for a review of recent papers).
The drawbacks of the actual schemes for diffusion are often linked with one or several of these key points:
\begin{itemize}
\item a non-local stencil (quite dense matrices);
\item cell-centered but also face-centered unknowns (large matrices);
\item non-symmetric definite positive matrices (loss of the energy balance);
\item loss of the convergence or of the accuracy on some particular grids;
\item loss of monotony for solutions in purely diffusive problems (the resulting matrix is not an ``M-matrix'').
\end{itemize}

We focus in this paper on the approximation of the Navier-Stokes and energy equations under the Boussinesq assumption, using a new scheme for diffusion terms.
This scheme is shown to provide a cell-centered approximation with a quite reduced stencil, leading to symmetric definite positive matrices and allowing a mathematical proof of convergence.
Although the diffusion matrix may not be shown to be an M-matrix in the general case, the maximum principle is nevertheless preserved in our numerical three-dimensional simulations. 
In this scheme, the discrete pressure gradient and the non-linear contributions are approximated so that 
the discrete kinetic and energy balances mimic their continuous counterparts. 
Indeed, the pressure gradient is chosen as the dual operator of the discrete divergence, and the discretization is such that there is no contribution of the non-linear velocity transport in the increase of kinetic energy.
In order to suppress the pressure checkerboard modes, the mass balance is stabilized by a pressure term which only redistributes the fluid mass within subsets of control volumes, the characteristic size of which is two or three times the local mesh size.\\
The remainder of this paper is divided into four sections. 
In Section 2, the continuous formulation is presented in the framework of free convection.
Section 3 presents the discrete scheme for general non-conforming meshes, with illustrations in the simplified case of  uniform rectangular grids,
and some mathematical properties.
The fourth section is devoted to the numerical validation, first with analytical solutions and 
then with  a classical natural convection problem.

\section{Continuous formulation}

Let $d$ be the dimension of the space ($d=2$ or $3$) and let $\O\subset\R^d$ be an open polygonal connected domain.
For $\bx\in\O$, our aim is to compute an approximation of the velocity $\bu(\bx)=\sum_{i=1}^du^{(i)}(\bx)\ei$, the pressure $p(\bx)$ and the temperature $T(\bx)$, solution of the steady and dimensionless Navier-Stokes and energy equations under the Boussinesq approximation:
\begin{subequations}\label{defNVSTKS}
\be \label{defNVSTKSa}
-\Pr \lap \bu + \grad p + (\bu\cdot\grad) \bu - \Ra\,\Pr\,T \ez=\bff(\bx) \mbox{ in } \Omega,
\ee
\be\label{defNVSTKSc}
- \lap T + (\bu\cdot\grad) T =g(\bx) \mbox{ in } \Omega,
\ee
\be\label{defNVSTKSb}
 \div \bu = 0\mbox{ in } \Omega,
\ee
\end{subequations}
where $\ez$ indicates the vertical upward direction, $\bff(\bx)=\sum_{i=1}^df^{(i)}(\bx)\ei$ and $g(\bx)$ are dimensionless regular functions modeling source or sink in the momentum or heat balances; $\Pr$ and $\Ra$  denote the Prandtl and Rayleigh numbers respectively.
We consider the case of the homogeneous Dirichlet boundary conditions for the velocity and
of the mixed Dirichlet-Neumann boundary conditions for the temperature. 
These boundary conditions read as follows:
\begin{equation}\label{diff_cls}
\left\{
\begin{array}{ll} 
\displaystyle \bu(\bx) =\b0 \qquad & \bx \in \Gamma, \\
\displaystyle T(\bx) =T_b(\bx) \qquad & \bx \in \Gamma_1, \\
\displaystyle -\grad T(\bx)\cdot\bn(\bx) =q_b(\bx) \qquad & \bx \in \frontomegadeux,
\end{array}\right.
\end{equation}
where $\frontomegaun, \frontomegadeux$ are subsets of the boundary $\frontomega$ of the domain $\Omega$ such that $\frontomegaun \cap \frontomegadeux=\emptyset $ and $\frontomegaun \cup \frontomegadeux=\frontomega$, and $\bn(\bx)$ is the outward unit normal vector to the boundary.
We assume that $T_b$ is the trace on $\frontomegaun$ of a function again denoted $T_b$ such that $T_b\in H^1(\O)$, and
we define the functional space $H_{\frontomegaun,0}^1(\O)=\{ T\in H^1(\O); T(\bx) =0 \mbox{ on } \frontomegaun\}$.
Then a weak formulation of equations (\ref{defNVSTKSa}-\ref{defNVSTKSb}) with boundary conditions (\ref{diff_cls}) reads: 

Find $\bu \in H_0^1(\O)^d$, $p\in L^2(\O)$ with $\int_\O p(\bx) \d \bx = 0$, 
and $T$ with $T-T_b \in H_{\frontomegaun,0}^1(\O)$ such that
\begin{subequations}\label{nstocontfsszero}
\be\label{nvstks_qqtmv}\ba
\dsp \Pr &\dsp \int_\O \grad \bu :\grad \bv \d \bx \dsp - \int_\O p \,\div \bv \d \bx 
+\int_\O \div(\bu\otimes\bu)\cdot\bv \d \bx\\
&\dsp - \Ra\,\Pr\int_\O \,T \ez \cdot \bv \d \bx=\int_\O \bff(\bx)\cdot \bv \d \bx,\ \forall \bv \in H^1_0(\O)^d,
\ea\ee
\be\ba
\dsp \int_\O \grad T \cdot \grad \theta \d \bx 
&+&\dsp \int_\O \div(\bu T)\theta \d \bx\\
 &=&\dsp \int_{\O} g(\bx) \theta \d \bx   
 -\int_{\frontomegadeux} q_b(\bx)  \theta(\bx) \d \bx, 
\ \forall \theta \in H_{\frontomegaun,0}^1(\O).\label{nvstks_nrj}
\ea\ee
\be
\div \bu(\bx) = 0 \hbox{ for a.e. } \bx\in\Omega,\label{nvstks_mass}
\ee
\end{subequations}
Although our discretization scheme belongs to the finite volume family, we shall also be using the weak form (\ref{nvstks_qqtmv}-\ref{nvstks_mass}) in our discretization.
Indeed, the discretization of the diffusive terms $-\Pr \lap \bu $ in \eqref{defNVSTKSa} and  $-\lap T $ in \eqref{defNVSTKSc} is obtained by the construction of a discrete gradient which is then used to approximate the term $ \Pr \int_\O \grad \bu :\grad \bv \d \bx  $ in \eqref{nvstks_qqtmv} and $\int_\O \grad T \cdot \grad \theta \d \bx$ in \eqref{nvstks_nrj}.

\section{Numerical scheme}

In this section we present the discretization scheme for Problem \eqref{defNVSTKS}-\eqref{diff_cls} under its weak form \eqref{nstocontfsszero}.  
The next paragraph is devoted to the notations for general discretization meshes and to the description of the discrete degrees of freedom.
We then describe the approximation of the  diffusive terms
(Sec. \ref{sec-vis}).
Because of the  collocated choice of the unknowns, a stabilization is needed. 
The stabilization we choose is imposed on the mass flux (rather than the overall balance) and also appears in   
 the momentum and energy equations through the convective contributions: this is described in Section \ref{sec-conv}.
It also involves the choice of some coefficients which are defined in Section \ref{sec-lambda}.  
The complete discrete problem is finally given in Section \ref{sec-complete}, and some of its mathematical properties sketched in Section \ref{sec-math}.

\subsection{Mesh and discrete spaces}

We denote by $\disc = (\mesh,\edges,\points)$ a space discretization, where (see Fig.~\ref{figmesh}):
\begin{itemize}
\item  $\mesh$ is a finite family of ``control volumes'', {\it i.e.} non empty connected open disjoint subsets of $\O$  such that $\overline{\O}= \displaystyle{\cup_{K \in \mesh} \overline{K}}$.
For any control volume $K\in\mesh$, we denote by $\dr K = \overline{K}\setminus K$  its boundary, $\mcv>0$  its measure (area if $d=2$, volume if $d=3$) and $h_K$ its diameter (that is the largest distance between any two points of $K$).
\item  $\edges$ is a finite family ``edges'' ($d=2$) or ``faces" ($d=3$) of the mesh; these are assumed to be non empty open disjoint subsets of $\overline{\O}$, which are included in a straight line if $d=2$  or in a plane if $d=3$, and with non zero measure. 
We assume that, for all $K \in \mesh$, there exists a subset $\edgescv$ of $\edges$ such that $\dr K = \displaystyle{\cup_{\edge \in \edgescv}}\overline{\edge} $. 
 The set $\edges$ is assumed to be partitioned into external and interior edges ($d=2$) or faces ($d=3$): $\edges=\edgesint\cup\edgesext$, with $  \edge   \subset \partial \Omega$ for any $ \edge \in \edgesext   $  and $ \edge   \subset  \Omega \setminus \partial \Omega $  for any $ \edge \in \edgesint$.
Any boundary edge $\edge$ is assumed to belong to a set $\edgescv$ for one and only one $K \in \mesh$;
any interior edge $\edge$ is assumed to belong to exactly two sets  $\edgescv$ and  $\edgescvv$ with $\cv \not = \cvv$, and in this case $\edge$ is included in the common boundary of $K$ and $L$, denoted  $K/L$. 
Note that there are cases in which $K/L$ includes two or more edges or faces, see for instance the third mesh for the unit cube, section \ref{sec-analyt}, and Figure \ref{fig:meshes}b.
We also assume that, if $\edge\in\edgesext$, then either $\edge\subset\frontomegaun$ or $\edge\subset\frontomegadeux$.
For all $\edge\in\edges$, we denote by $\xedge$ and $\medge$ the barycenter and the measure of $\edge$.
For all $K \in \mesh$ and $\edge \in \edgescv$, we denote by $\ncvedge$ the unit vector normal to $\edge$ outward to $K$.
\item $\points$ is a family of collocation points $\points = (\xcv)_{K \in \mesh}$ of $\O$ which is chosen such that for all $K\in\mesh$ and for all $\bx\in K$, the property $[\xcv,\bx]\subset K$ holds.
Note that this choice is possible for quite general polygons, including those with re-entrant corners, see Fig. \ref{figmesh}.
The Euclidean distance $d_{K,\edge}$ between $\xcv$ and the hyperplane including $\edge$ is thus  positive.
We also denote by $C_{K,\edge}$ the cone with vertex $\xcv$ and basis $\edge$.
\end{itemize}
 Next, for any $\edge \in \edgesint$, we choose some real coefficients $(\beta_\edge^\cvv)_{\cvv\in\mesh}$ such that  the barycenter    $\xedge$ of $\edge$ is expressed by
\be
\ \xedge = \sum_{\cvv\in\mesh} \beta_\edge^\cvv\xcvv, \qquad \sum_{\cvv\in\mesh} \beta_\edge^\cvv = 1.
\label{bary}\ee
In three space dimensions, it is always possible to restrict the number of nonzero coefficients $\beta_\edge^\cvv$ to four (in practice, the scheme has been shown to be robust with respect to the choice of these four control volumes, taken close enough to the considered edge).
\begin{itemize}
\item[] {\em Note that in the case of an uniform rectangular grid which is depicted in Figure \ref{figmeshcart}, an obvious choice for the coefficients $\beta_\edge^\cvv$ is obtained by noticing that ${\bf x}_{i+1/2,j} = ({\bf x}_{i,j} + {\bf x}_{i+1,j})/2$ and ${\bf x}_{i,j+1/2} = ({\bf x}_{i,j} + {\bf x}_{i,j+1})/2.$ Thus for any edge $\sigma$, only two coefficients $\beta_\edge^\cvv$ need to be nonzero.}
\end{itemize}

We now define the finite dimensional space $\R^\mesh\times \R^\edges$ (where an element $\tstf\in\R^\mesh\times \R^\edges$ is defined by the family of real values $((\tstf_\cv)_{\cv\in\mesh},(\tstf_\edge)_{\edge\in\edges})$) and the following subspaces:
\begin{itemize}
\item $\dset=\left\{u\in \R^\mesh\times \R^\edges,\ \forall \edge\in\edgesint, u_\edge= \sum_{\cvv\in\mesh} \beta_\edge^\cvv  u_\cvv \right\}$ 
(the dimension of $\dset$ is the number of control volumes plus that of boundary edges),
\item $\dsetU=\left\{u\in \dset,\ \forall \edge\in\edgesext, u_\edge=0 \right\}$
(the dimension of $\dsetU$ is the number of control volumes),
\item $\dsetTheta=\left\{\theta\in \dset,\ \forall \edge\in\edgesext \cap\frontomegaun, \theta_\edge=0 \right\}$
(the dimension of $\dsetTheta$ is the number of control volumes plus that of boundary edges on $\frontomegadeux$).
\end{itemize}

\subsection{Discretization of  diffusive terms}\label{sec-vis}

Let us first define a discrete gradient for the elements of $\dset$ on cell $\cv\in\mesh$.
We set, for any $u\in\dset$ and $\cv\in\mesh$:
\begin{equation}
\label{defgradcv}
\grad_\cv u =\dsp \frac 1 \mcv \dsp \sum_{\edge\in\edgescv} \medge (u_\sigma - u_\cv) \ncvedge.
\end{equation}
Note that this is a centered gradient.
\begin{itemize} 
	\item[] \em As an illustration, consider the case of the two dimensional uniform rectangular grid depicted in Fig. \ref{figmeshcart}, let $u \in\dset$, and choose the natural choice $\beta_\sigma^L = 1/2$  if $\sigma$ is a side of $L$ and $0$ otherwise. 
Then with the (natural) notations of Fig. \ref{figmeshcart}, one has:
$$
\grad_{\cv_{i,j}} u = \left(
\begin{array}{c} \dsp \frac 1 {2h_x} (u_{i+1,j} - u_{i-1,j}) \\ \dsp \frac 1 {2h_y} (u_{i,j+1} - u_{i,j-1})\end{array}\right);
$$
if we apply this formula to the  element $1_{ \cv_{i,j} } \in\dset$ with all of  components equal to zero except for the one associated with $\cv_{i,j}$ which is equal to 1, we get that
the vector $\nabla_\cvv { 1}_{\cv_{i,j}}$ is zero for all control volumes $\cvv$ except those neighboring ${\cv_{i,j}}$, as shown on  Fig. \ref{figgradcart}a. 
Considering the checkerboard solutions on uniform rectangular grids,
the above expression shows that this discrete gradient may vanish for some non-constant functions. 
\end{itemize}
An approximation of the diffusive terms $ \Pr \int_\O \grad \bu :\grad \bv \d \bx  $ in \eqref{nvstks_qqtmv} and $\int_\O \grad T \cdot \grad \theta \d \bx$ in \eqref{nvstks_nrj} using this discrete gradient (\ref{defgradcv}) would yield a non-coercive form.
Thus, we shall work with a modified gradient, defined in (\ref{grad2})-\eqref{defgrad} below.

To this end,  for all $\edge\in\edgescv$ we first define $R_{\cv,\edge} u\in\R$ which may be seen as a consistency error on the normal flux, by:
\[
R_{K,\edge} u = \frac {\sqrt{d}} {\dcvedge}  \left( u_\edge - u_K - \grad_K u \cdot(\xedge-\xcv)\right).
\label{defresidu}
\]
One may note that $R_{K,\edge} u=0$ if $u_K$ and $u_\edge$ are the exact values of a linear function at points $\xcv$ and $\xedge$, for all $\cv$ and $\edge\in\edgescv$.
We then give the following expression for a stabilized discrete gradient of $u\in\dset$ in each cone $C_{\cv,\edge}$:
\be\label{grad2}
\grad_{\cv,\edge} u =\grad_\cv u + R_{K,\edge} u\,\ncvedge. 
\ee
\begin{itemize}
	\item[] \em In the case of the uniform rectangular grid given in Fig. \ref{figmeshcart} and for $K = K_{i,j}$ and $\sigma = \sigma_{i+1/2,j}$, we find:
$$
R_{K_{i,j},\edge_{i+1/2,j}}= \dfrac{\sqrt{2}}{2h_x}(u_{i+1,j}+u_{i-1,j}-2u_{i,j}).
$$
The stabilized discrete gradient (\ref{grad2}) applied  to the above defined element $1_{\cv_{i,j}}$ of $\dset$ provides nonzero contributions on the triangular subcells $C_{K,\sigma}$ of the cell $\cv_{i,j}$ and of its neighbors (Fig. \ref{figgradcart}b).
\end{itemize}
The global discrete gradient is chosen as the function $\grad_\disc u$:
\begin{equation}
\label{defgrad}
\grad_\disc u(\bx) =\grad_{\cv,\edge} u,\mbox{ for a.e. }\bx\in C_{\cv,\edge},\ \forall \cv\in\mesh,\ \forall \edge\in\edgescv.
\end{equation}
 We then plan to approximate the term $\int_\Omega \grad  u(\bx) \cdot \grad  v(\bx) \d \bx$  
 by:
\be
\int_\O \grad_\disc u(\bx) \cdot \grad_\disc v(\bx) \d \bx =
\sum_{\cv\in\mesh} \sum_{\edge\in\edgescv}\frac{\medge\dcvedge}{d} \grad_{\cv,\edge} u \cdot \grad_{\cv,\edge} v,\ \forall u,v\in\dset.
\label{defbilinear}\ee
In fact, it is shown in  \cite{crasEH,eym-08-dis} that this expression defines a symmetric inner product on $\dset$, and provides a good approximation for $\int_\O \grad u(\bx) \cdot \grad v(\bx) \d \bx$; this approximation may be seen as a low degree discontinuous Galerkin method.
If one seeks a finite volume interpretation of this scheme, it is possible, expressing $u_\edge$ and $v_\edge$ for all $\edge\in\edgesint$ thanks to the 
relations (\ref{bary}), to show that
\be\ba
\dsp\int_\O \grad_\disc u(\bx) \cdot \grad_\disc v(\bx) \d \bx = \\
\dsp\sum_{K\in\mesh}\left( \sum_{L\in {\cal N}_K} F_{\cv,\cvv}(u) v_\cv + \sum_{\edge\in\edgescv\cap \edgesext} F_{\cv,\edge}(u) ( v_K - v_\edge)\right),
\ea\label{propbilinear}\ee
 where for any $K\in\mesh$,  ${\cal N}_\cv$ is the subset of cells playing a part in the barycenter expression of $\xedge$, for all edges of the cell $\cv$ and of the neighbors of $\cv$, {\it i.e.} ${\cal N}_\cv= \{\cvvv\in\mesh;\beta_\edgep^\cvvv \neq 0,\forall\edgep\in\edgescvv,\forall\cvv\in{\cal M}_\edge,\forall\edge\in \edgescv\}$;
$F_{\cv,\cvv}(u)$ is a linear function of the unknowns $(u_L)_{L \in \mesh}$  which is such that $F_{\cv,\cvv}(u) = - F_{\cvv,\cv}(u)$.
In the general case, the expression of $F_{\cv,\cvv}(u)$ is rather complicated (see \cite{eym-08-dis}).
\begin{itemize}
	\item[] \em In the case of the uniform rectangular grid of Fig. \ref{figgradcart}, this expression simplifies into the usual two point flux; for instance, the flux from $\cv_{i,j}$ to $\cv_{i+1,j}$ reads:
$$
F_{\cv_{i,j},\cv_{i+1,j}}(u)   = h_y \dfrac{u_{i,j} - u_{i+1,j}}{h_x}.
$$
More generally, a two point flux is also  obtained in the two or three dimensional non uniform rectangular cases. 
Indeed, locating $\xcv$ at the center of gravity of the cell $\cv$, the relation $(\xedge-\xcv)/\dcvedge=\ncvedge$ holds. 
It is then possible to write ${\bf x}_\edge =(\dcvvedge {\bf x}_\cv+\dcvedge {\bf x}_\cvv)/(\dcvvedge+\dcvedge)$ and 
$u_\edge(\tstf)=(\dcvvedge u_\cv+\dcvedge u_\cvv)/(\dcvvedge+\dcvedge)$ for all $\edge$ such that $\edge \subset K/L$, and for all $u\in\dset$. 
Using the identity 
\[\label{rel_magic}
\sum_{\edge\in\edgescv}\medge(\bx_\edge-\bx_\cv)\ncvedge^t=\mcv \mathbb{I}
\]
where $^t$ designates the transposition and $\mathbb{I}$ the identity matrix, we obtain \cite[Lemma 2.1] {eym-08-dis}:
\[
\ba
\dsp 
\int_\O \grad_\disc u(\bx) \cdot \grad_\disc v(\bx) \d \bx =&
\dsp \sum_{\edge\in\edgesint,\edge \subset K/L} \frac{\medge}{\dcvedge+\dcvvedge}(u_\cvv-u_\cv)(v_\cvv-v_\cv)\\ 
&+\dsp \sum_{\edge\in\edgesext \cap\edgescv} \frac {\medge} {\dcvedge}(u_\edge - u_{K}) (v_\edge - v_{K}).
\ea
\]
Then the previous relation leads to define ${\cal N}_K$ as simply the set of the natural neighbors of $K$, and to define the fluxes by the natural two-point difference scheme, in the same manner as in \cite{Chenier05,Touazi2007}:
\begin{equation}
\label{diffussion_scheme_superadm}\ba
F_{\cv,\cvv}(u) &\dsp = \frac{\medge}{\dcvedge+\dcvvedge}(u_\cv-u_\cvv)\ & \mbox{ for } \edge\in\edgesint,   \edge \subset K/L\\
 F_{\cv,\edge}(u) &\dsp= \frac {\medge} {\dcvedge}(u_K - u_\edge)\ & \mbox{ for } {\edge\in\edgesext\cap\edgescv.}
\ea
\end{equation}
The classical and cheap 5- and 7-point schemes on rectangular or orthogonal parallelepipedic meshes is then recovered.
An advantage can then be taken from this property, by using meshes which consist in orthogonal parallelepipedic control volumes in the main part of the interior of the domain, as illustrated by the cone-shaped cavity (Fig. \ref{fig:meshes}c) in Section \ref{sec-analyt}. 
\end{itemize}

Note that the approximation of $-\int_K \Delta u \d\bx$ is obtained by letting $v_K = 1$, $v_L = 0$ for $L\neq K$ and $v_\edge = 0$ for $\edge\in\edgesext$ in (\ref{propbilinear}):
\[
-\int_K \Delta u \d\bx \simeq \sum_{L\in {\cal N}_K} F_{\cv,\cvv}(u) + \sum_{\edge\in\edgescv\cap \edgesext} F_{\cv,\edge}(u),
\label{bilanlap}
\]
so that we may define an approximate Laplace operator $ \Delta_\disc $ by the constant values $ \Delta_K u$ on the cells $K$:
\be
- \Delta_K u = \frac 1 \mK \left( \sum_{L\in {\cal N}_K} F_{\cv,\cvv}(u) + \sum_{\edge\in\edgescv\cap \edgesext} F_{\cv,\edge}(u)\right).
\label{lapdis}
\ee
\begin{itemize}
	\item[] \em In the case of the uniform rectangular grid of Fig. \ref{figgradcart}, this discrete Laplacian leads to the usual five point formula: 
$$
- \Delta_{\cv_{i,j}} u = \frac 1 { {h_{x}}^2} \left( 2 u_{i,j} -  u_{i+1,j} - u_{i-1,j}\right) +  \frac 1 { {h_{y}}^2} \left( 2 u_{i,j} -  u_{i,j+1} - u_{i,j-1}\right).
$$
\end{itemize}
In the general case, the stencil of the discrete operator on cell $\cv$ is defined by ${\cal N}_\cv$ (see relation (\ref{lapdis})) and therefore depends on the way the barycenters $\xedge$ are computed.
For general grids, the equation for a given cell usually concerns the unknowns associated to itself, its neighbors, the neighbors of its neighbors and possibly some additional adjacent cells.
The resulting matrix is usually not an ``M-matrix'', except on particular meshes such as conforming orthogonal parallelepipeds, in which case we obtain the usual two point flux scheme, as previously pointed out.

\subsection{Pressure-velocity coupling, mass balance and convective contributions}\label{sec-conv}
For all $\bu\in (\dsetU)^d$, we define a discrete divergence operator 
 by:
\be \label{eq:divergence}
\div_K \bu = \frac{1}{\mcv}\sum_{\edge\in\edgescv} \medge \bu_\edge\cdot \ncvedge,\qquad \forall \cv\in\mesh.
\ee
where $\bu_\edge=\sum_{i=1}^d u_\edge^{(i)} \ei$.
Notice that  
\[
\div_K \bu = \sum_{i=1}^d (\grad_K \ui(i))^{(i)}, 
\]
with $\grad_K \ui(i)$ defined by \eqref{defgradcv}.
\begin{itemize}
	\item[] \em In the case of the uniform rectangular mesh of Fig. \ref{figmeshcart}, this operator reads:
$$
\div_{K_{i,j}} \bu = \dfrac{1}{2h_x} \left(u_{i+1,j}^{(1)} - u_{i-1,j}^{(1)}\right) + \dfrac{1}{2h_y} \left(u_{i,j+1}^{(2)} - u_{i,j-1}^{(2)}\right). $$
\end{itemize}
We then define the function $\div_\disc \bu$ by the relation
\[
\div_\disc \bu (\bx) = \div_K \bu, \mbox{ for a.e. }\bx\in K,\ \forall K\in\mesh.
\]

The discrete gradient operator used for the pressure gradient is defined as the dual of this divergence operator.
More precisely, we mimic at the discrete level the (formal) equality $\int_\Omega p \ \div \bv \ \d \bx = - \int_\Omega \grad p \cdot \bv \ \d \bx$.  
The discrete equivalent of $\int_\Omega p \ \div \bv  \ \d \bx$ reads $\sum_{\cvv \in \mesh} m_L \ p_L \ \div_\cvv \bv $ with $\bv \in (\dset)^d $ and $p \in \dset$; we then define the discrete pressure gradient which we denote  $\widehat\grad_\cv p$ (on cell $\cv$) , such that 
$$\sum_{\cvv \in \mesh} m_\cvv \widehat\grad_\cvv p \cdot \bv_\cvv  = -\sum_{\cvv\in\mesh} \mcvv p_\cvv  \div_\cvv \bv , \forall \ \bv \in (\dset)^d.$$
>From the definition of the divergence \eqref{eq:divergence}  and  of $\dset$, we thus seek  $(\widehat\grad_\cvv p)_{\cvv \in \mesh}$ such that:
\be \label{eq:gradpweak}
\sum_{\cvv \in \mesh} m_L \widehat\grad_\cvv p \cdot \bv_\cvv  = -\sum_{\cvv\in\mesh}p_\cvv \sum_{\edge\in\edgescvv\cap\edgesint}\medge \sum_{\cvvv\in\mesh} \beta_\edge^\cvvv \bv_\cvvv \cdot \ncvvedge.
\ee
Taking for $\bv$ the element of $(\dset)^d$ with components $v_\cvv^{(j)}=1$ if $j=i$ and $\cvv=\cv$, and $0$ otherwise , we thus get:
\begin{align}\label{eq:gradp}
\mcv \widehat\grad_\cv p &=  -\sum_{\cvv\in\mesh}p_\cvv \sum_{\edge\in\edgescvv\cap\edgesint}\medge  \beta_\edge^\cv  \ncvvedge \nonumber \\
&=
\sum_{\edge\in\edgesint,\edge \subset \cvv/\cvvv} \medge  \beta_\edge^\cv (p_\cvvv-p_\cvv)\ncvvedge.
\end{align}
Remark that  $\widehat\grad_\cv p $  is neither constructed with the discrete gradient (\ref{defgradcv}) nor with the stabilized one (\ref{grad2});
its expression is the dual form of the divergence (\ref{eq:divergence}).\\
\begin{itemize}
	\item[] \em 
	In the case of non uniform rectangles ($d=2$) or parallelepipeds ($d=3$) with collocated points at the gravity center of the cells, and for all $\edge\in\edgesint$ with $\edge\subset\cv/\cvv$, we only need two  non-zero coefficients $\beta_\edge^\cvvv$: 
$\beta_\edge^\cv=\frac{\dcvvedge}{\dcvedge+\dcvvedge}$ and $\beta_\edge^\cvv=\frac{\dcvedge}{\dcvedge+\dcvvedge}$.
Therefore, relation (\ref{eq:gradp}) reduces to 
$$
\mcv \widehat\grad_\cv p =  
\sum_{\edge\in\edgescv,\edge \subset \cv/\cvv} \medge  \frac{\dcvvedge}{\dcvedge+\dcvvedge} (p_\cvv-p_\cv)\ncvedge.
$$
Using  $p_\cv\sum_{\edge\in\edgescv}\medge\ncvedge=0$,
$$
\mcv \widehat\grad_\cv p = 
\sum_{\edge\in\edgescv,\edge \subset \cv/\cvv} 
\medge\frac{\dcvvedge p_\cvv+\dcvedge p_\cv}{\dcvedge+\dcvvedge}\ncvedge  +
\sum_{\edge\in\edgescv\cap\edgesext} \medge p_\cv\ncvedge
$$
	For a uniform grid and a control volume without boundary faces, the above expression resumes to
$$
\mcv\widehat \grad_\cv p = \sum_{\edge\in\edgescv, \edge \subset K/L} 
\medge\frac{p_\cvv+p_\cv}{2}\ncvedge,
$$
which provides, in the particular case of Figure \ref{figmeshcart}, the usual formulation
$$
h_xh_y \widehat\grad_\cv p = h_y\frac{p_{i+1,j}-p_{i-1,j}}{2}+h_x\frac{p_{i,j+1}-p_{i,j-1}}{2}.
$$
\end{itemize}

As recalled in the introduction of this paper, a pressure stabilization method is implemented in the mass conservation equation in order
to prevent from oscillations of the pressure, as for instance in \cite{bre-84-sta} in the finite element setting, \cite{Nagele2007,Rhie1983} in the finite volume setting. 
The originality of our approach is that we directly include the stabilizing diffusive pressure flux in the approximated mass flux, so that it will appear not only (as usual) in the mass equation, but also in the momentum equation through the non-linear convective term. From the mathematical point of view, this helps in obtaining  simple estimates on the velocity and pressure, but more importantly, it ensures that the contribution of the discrete non-linear convective term to the kinetic (and thermal) energy balance is zero, just as in the continuous case. 
Let us define the stabilized mass flux across  $\edge\subset K/L$ by
\be\label{phi}
\Phi_{K,\edge}^{\lambda}(\bu,p) = \medge \left(\bu_\edge\cdot \ncvedge + \lambda_\edge (p_K- p_L)\right),
\ee
where $(\lambda_\edge)_{\edge\in\edgesint}$ is a given family of positive real numbers, the choice of which is discussed below.
Note that the quantity $\lambda_\edge (p_K- p_L)$ may be seen as a numerical pressure diffusion flux, and that the overall numerical flux remains conservative, that is, if    $\edge\subset K/L$  then $\Phi_{K,\edge}^{\lambda}(\bu,p) + \Phi_{L,\edge}^{\lambda}(\bu,p) = 0$.
\begin{itemize}
	\item[]\em 
	In the case of the uniform rectangular mesh of Fig. \ref{figmeshcart}, the expression of this flux through a vertical edge $\edge_{i+1/2,j}$ reads:
$$
\Phi_{K_{i,j},\edge_{i+1/2,j}}^{\lambda} \bu = \dfrac{h_y}{2} \left(u_{i+1,j}^{(1)} + u_{i,j}^{(1)}\right)   + h_y \lambda_{i+1/2,j}  (p_{i,j} - p_{i+1,j}). 
$$
\end{itemize}
We then use the modified flux, in order to define a stabilized centered transport operator which is defined, for all $\bu\in (\dsetU)^d$, $w\in\dset$ and $\cv\in\mesh$, by
\[
\div_K^{\lambda}(w,\bu,p) = \frac{1}{\mcv}
\sum_{\edge\in\edgescv, \edge \subset K/L} \Phi_{K,\edge}^{\lambda}(\bu,p) \frac {w_K + w_L} 2.
\]
An interesting remark is that, in the case where the mass balance equation in the control volume $K$ is satisfied, that is:
\[ \label{discret_mass}
\div_K^{\lambda}(1,\bu,p) = \frac 1 {\mcv} 
\sum_{\edge\in\edgescv\cap\edgesint}\Phi_{K,\edge}^{\lambda}(\bu,p)=0,
\]
 then 
\[ \label{anti-bernoulli}
\sum_{\edge\in\edgescv\cap\edgesint}\Phi_{K,\edge}^{\lambda}(\bu,p) w_K=0,
\]
so that the following relation also holds: 
\[ \label{transport}\ba
\dsp  \div_K^{\lambda}(w,\bu,p) = \frac 1 {\mcv}\sum_{\edge\in\edgescv, \edge \subset K/L} \Phi_{K,\edge}^{\lambda}(\bu,p) \frac {w_L - w_K} 2.
\ea\] 
We shall use this latter form in the practical implementation, in particular in the discretization of the non-linear convection term. 
Indeed, it is more efficient when computing the Jacobian matrix of the  momentum equation, since it avoids summing up values of the same amplitude.
\begin{itemize}
	\item[]\em
	In the particular case of the uniform rectangular grid of Fig. \ref{figmeshcart}, the summation $\sum_{\edge\in\edgescv, \edge \subset K/ L}$ involves the four edges  between the control volume $K = K_{i,j}$ and its neighbors $L = K_{i+1,j}$, $K_{i-1,j}$, $K_{i,j+1}$, $K_{i,j-1}$. 
	If $L = K_{i-1,j}$ the 
expression which appears in the summation reads:
$$ 
h_y \left(  - \dfrac{u_{i,j}^{(1)}+u_{i-1,j}^{(1)}}{2} + \lambda_{i-1/2,j}(p_{i,j}-p_{i-1,j}) \right)  
\dfrac{w_{i,j} +w_{i-1,j} }{2}.
$$ 
\end{itemize}

When the local grid Reynolds (or P\'eclet) number is much larger than $1$, an upwind scheme must be applied that consists in substituting $\div_K^{\lambda}(w,\bu,p)$ by 
\[\ba
\div_K^{\lambda,{\rm up}}(w,\bu,p) = \\ \dsp
 \frac{1}{\mcv}
\sum_{\edge\in\edgescv, \edge \subset K/L} \big(\max(\Phi_{K,\edge}^{\lambda}(\bu,p),0) w_K + 
\min(\Phi_{K,\edge}^{\lambda}(\bu,p),0) w_L\big).
\ea\]
In both the centered or upwind cases, the functions $\div_\disc^{\lambda}(w,\bu,p)$ and $\div_\disc^{\lambda,{\rm up}}(w,\bu,p)$ are defined by their constant values in each control volume.
For $\bu,\bw \in (\dsetU)^d$, we  also define the centered vector divergence operator $\bdiv_\disc^{\lambda}(\bw,\bu,p)$  such that the $i$-th component of  $\bdiv_\disc^{\lambda}(\bw,\bu,p)$ is equal to  $\div_\disc^{\lambda}(w^{(i)},\bu,p)$, for $i = 1,\ldots,d$;  a similar expression applies for the upwind vector divergence  operator $\bdiv_\disc^{\lambda,{\rm up}}(\bw,\bu,p)$.

\subsection{Choice for the parameters $(\lambda_\edge)_{\edge\in\edgesint}$ }\label{sec-lambda}
\begin{figure}[htbp] 
 \centering 
 \begin{minipage}{.45\textwidth}
 \includegraphics[width=\textwidth]{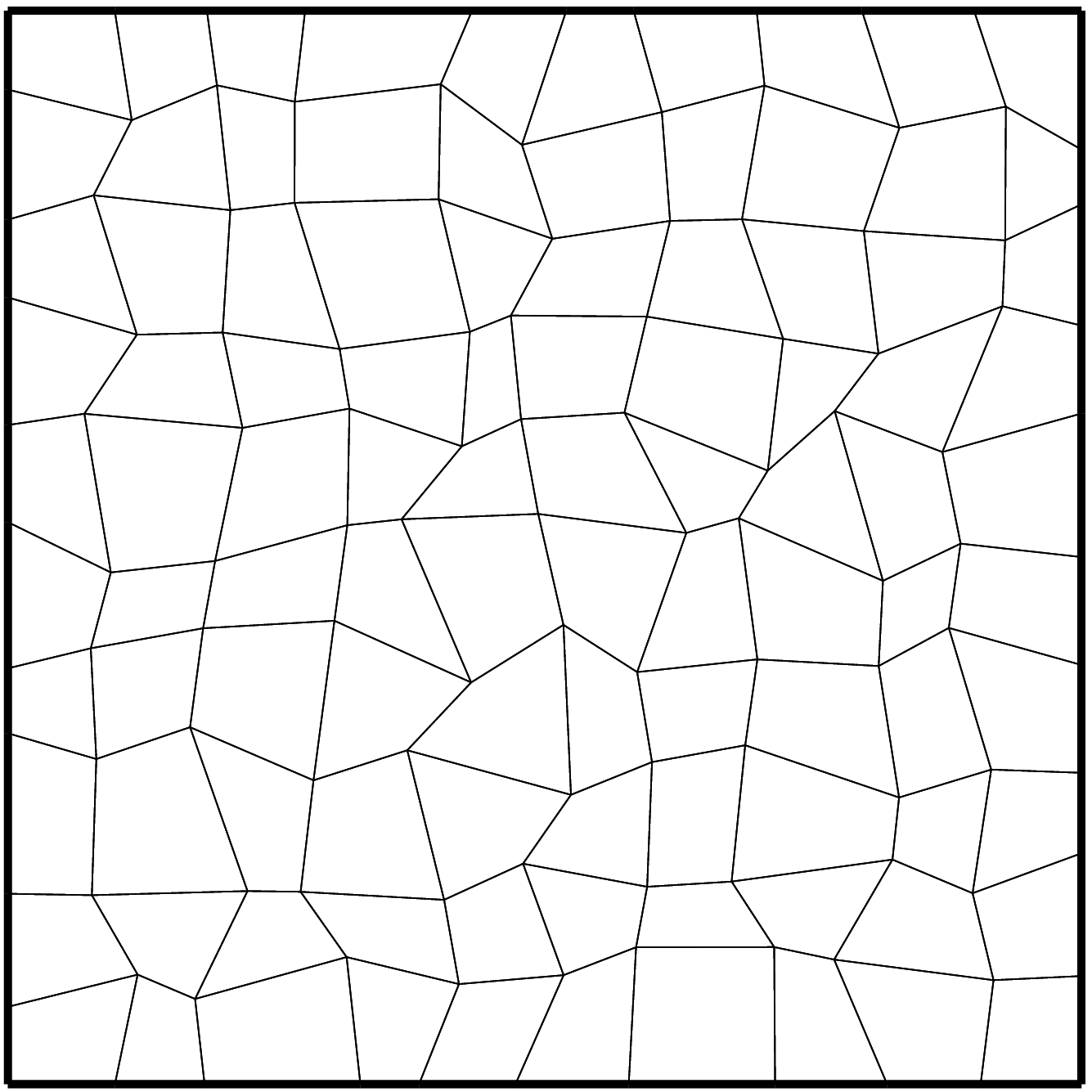}
 \centering (a) 
 \end{minipage}
 \begin{minipage}{.45\textwidth}
 \includegraphics[width=\textwidth]{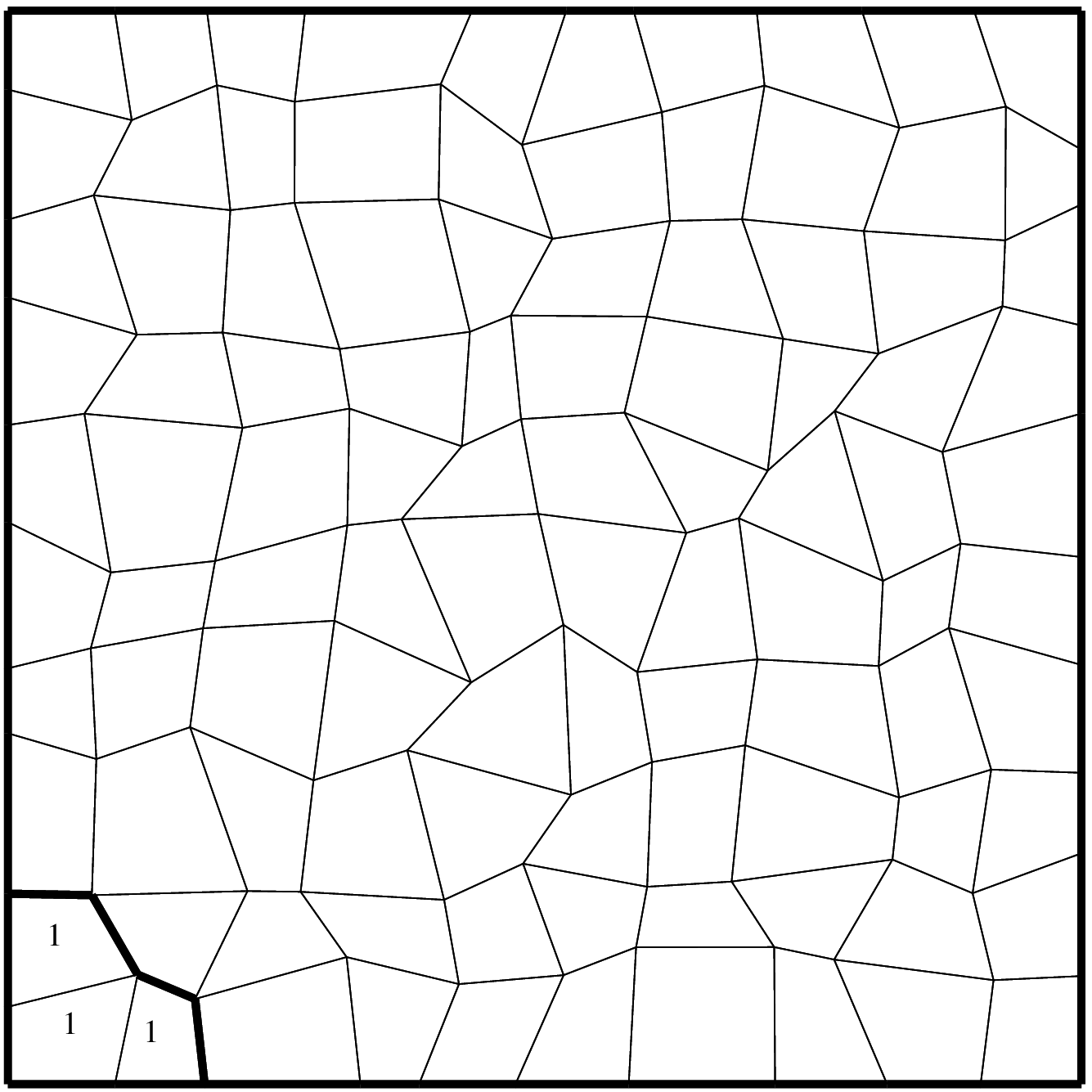}
 \centering (b) 
 \end{minipage}
 \begin{minipage}{.45\textwidth}
 \includegraphics[width=\textwidth]{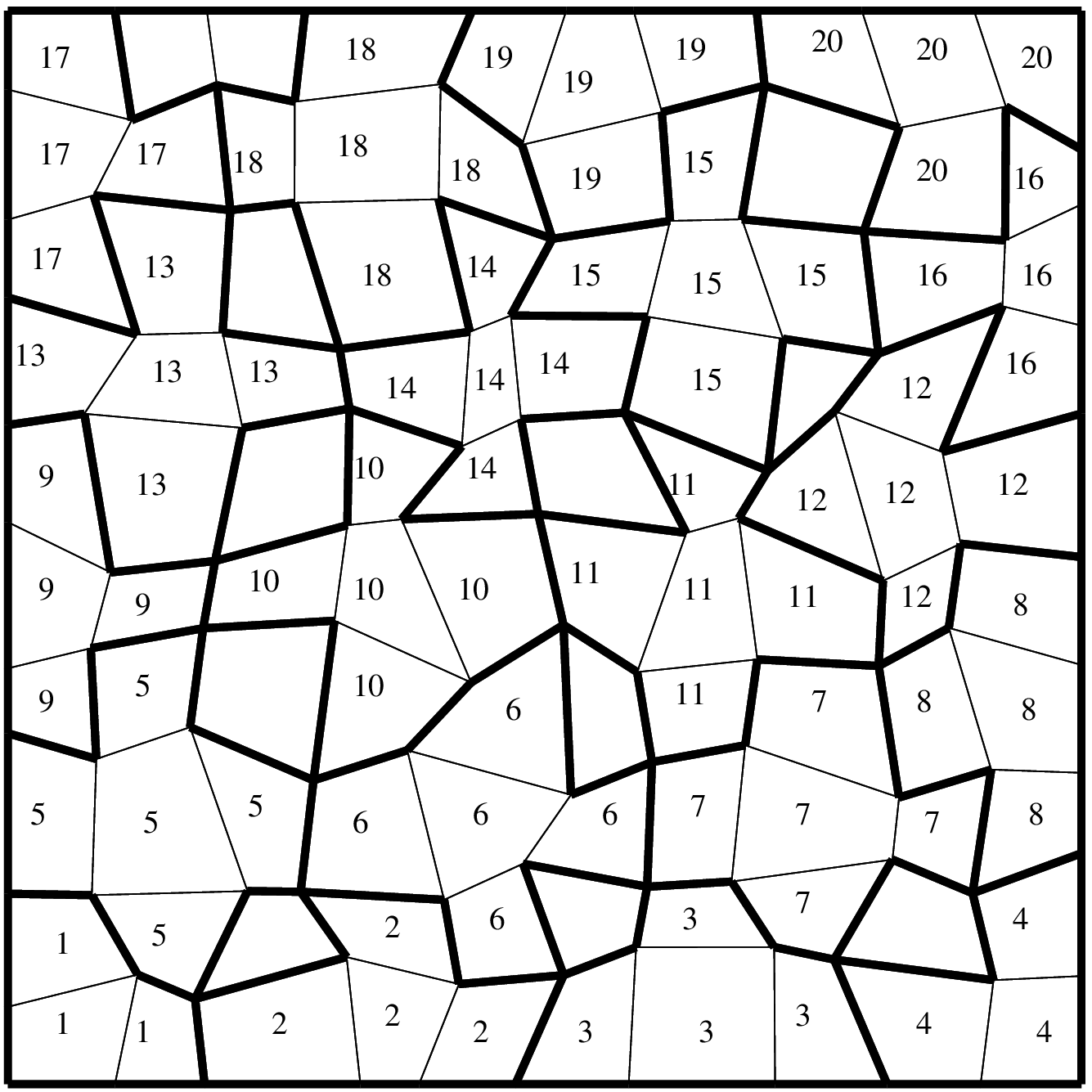}
 \centering (c) 
 \end{minipage} 
 	 \begin{minipage}{.45\textwidth}
 \includegraphics[width=\textwidth]{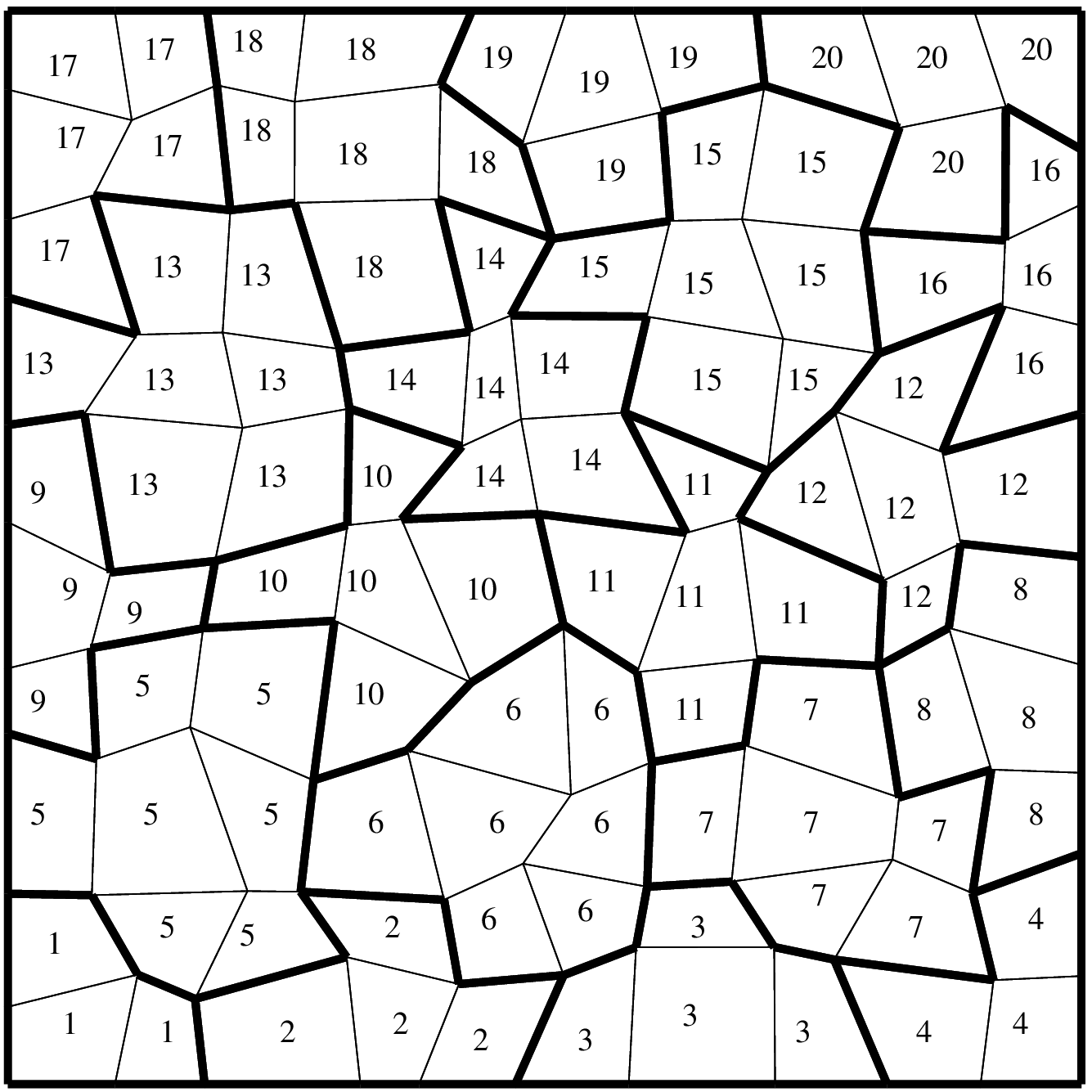}
 \centering (d) 
 \end{minipage}
 \caption[Example of cluster construction; (a) Mesh, (b) First cluster, (c) All the clusters are constructed but isolated cells (unnumbered cells) still remain, (d) The isolated cells are connected to  neighboring clusters.]{}
 \label{fig:cluster}
\end{figure}

Different strategies can be applied to define the parameters $(\lambda_\edge)_{\edge\in\edgesint}$.
Amongst all of them we applied the "cluster stabilization method" \cite{Chenier05,Touazi2007}; 
it consists in constructing a partition of $\mesh$, denoted $\clusters$, and setting 
$\lambda_\edge = \lambda>0$ if there exists $G\in\clusters$ (such $G\subset\mesh$ is called a cluster) with $\edge\subset K/L$,
$K$ and $L$ belonging to $G$, and $\lambda_\edge = 0$ otherwise.
{ Here is an example of an algorithm creating a cluster partition:
\begin{enumerate}
  \item for all cells $\cv\in\mesh$, initialize a new cluster if $\cv$ and its neighboring cells 
do not already belong to a cluster;
	\item for any remaining isolated cell $\cvv$, connect it to the closest cluster having the largest number of common edges with $\cvv$.
\end{enumerate} 
This algorithm is now applied to the mesh of figure (\ref{fig:cluster}a), the cells being described from left to right, from the lower to the upper row.
Figures (\ref{fig:cluster}b) and (\ref{fig:cluster}c) illustrate the first cluster  and the set of the clusters at the end of the first step of the algorithm.
Figure \ref{fig:cluster}d shows the clusters after the isolated (unnumbered) cells of figure \ref{fig:cluster}c have been connected.
The choice of the value $\lambda$ depends on the physical problem and it must be chosen both large to avoid the appearance of spurious modes and small to preserve an accurate approximation of the mass equation. 
In \cite{Touazi2007}, comparisons with other stabilization methods were performed and the results showed that solutions are little sensitive to the value of  $\lambda$.
For the natural convection example presented in section~\ref{Natural_convection_problem}, $\lambda=10^{-8}$.
}

\subsection{Resulting discrete equations}\label{sec-complete}

We denote by $T_{b,\disc}$ the element $T\in\dset$ such that $T_K=0$ for all $K\in\mesh$,
$T_\edge = 0$ for all $\edge\in\edgesint$ and all $\edge\in\edgesext$ with $\edge\subset\frontomegadeux$, and, for
 all $\edge\in\edgesext$ with $\edge\subset\frontomegaun$,
\be\label{nvstksd_nrjcld}\ba
\dsp T_\edge = \frac{1}{\medge}\int_\edge T_b(\bx) \d s(\bx).
\ea\ee
Let $\hmesh\subset L^2(\O)$ denote the set of functions which are constant in each $\cv\in \mesh$; for any function $q\in \hmesh$, we shall denote by $q_\cv$ its constant value on $\cv\in\mesh$.
We then define the mapping $P_{\mesh}~:~\dset\to \hmesh$ by $\tstf\in\dset \mapsto P_{\mesh} \tstf$ with  $ P_{\mesh} \tstf(\bx) = \tstf_{\cv}$ for a.e. $\bx\in \cv$ and all $\cv\in\mesh$.
We also define the mapping $P_{\edges}~:~\dset\to L^2(\frontomega)$ by
 $\tstf\in\dset\mapsto P_{\edges} \tstf $ with  $ P_{\edges} \tstf (\bx) = \tstf_{\edge}$ for a.e.  $\bx\in \edge$ and all $\edge\in\edgesext$.

Let us then use the previously defined discrete operators to formulate a discrete approximation to problem \eqref{nstocontfsszero}:  

Find $\bu = (\ui(i))_{i=1,d} \in (\dsetU)^d$, 
$p\in\hmesh$ with $\int_\O p(\bx)\d\bx=\sum_{\cv\in\mesh}\mcv p_\cv=0$ and $T - T_{b,\disc} \in \dsetTheta$ such that: 
\be\ba
\dsp \Pr\int_\O \grad_\disc \bu : \grad_\disc \bv \ \d \bx
-\int_\O p \ \div_\disc\bv \ \d\bx + \int_\O \bdiv_\disc^{\lambda}(\bu,\bu,p)\cdot  P_{\mesh} \bv  \ \d \bx\\
\dsp - \Ra\,\Pr \int_\O P_{\mesh} T \ \ez \cdot P_{\mesh} \bv \ \d \bx = \int_\O \bff\cdot P_{\mesh} \bv \ \d \bx,\
\forall \bv \in (\dsetU)^d,
\ea\label{eqintmomb}\ee
 
\be\ba
\dsp\int_\O \grad_\disc T \cdot \grad_\disc \theta \ \d \bx + \int_\O \div_\disc^{\lambda}(T,\bu,p) P_{\mesh} \theta \ \d \bx \\
\dsp = \int_\O g \ P_{\mesh} \theta \ \d \bx - \int_{\frontomegadeux} q_b P_{\edges} \theta \ \d s, 
 \ \forall \theta \in \dsetTheta,
\ea\label{eqenerb}\ee
 
\be\ba
\div_\disc^{\lambda}(1,\bu,p) = 0\mbox{ a.e. in }\O.
\ea\label{eqmasb}\ee

We then deduce from (\ref{eqintmomb}) the $d$ discrete  momentum balances over the control volume $K$, letting $\vi(i)=1$ in $K$, and $0$ otherwise; these equations read, in vector form:

\be\label{nvstksd_qqtmv}\ba
-\Pr\,\mK \Delta_K \bu 
\dsp +\sum_{\edge\in\edgesint,\edge \subset M/L } \medge \beta_\edge^K (p_M-p_L) \ \ncvvedge\\
\dsp + \sum_{\edge\in\edgescv, \edge \subset K/L} \Phi_{K,\edge}^{\lambda} (\bu,p) \frac{\bu_K+\bu_L} 2 \dsp - \Ra\,\Pr\ \mcv T_K \ez = \int_K\bff \ \d \bx
\ea\ee

(where $-\Delta_K \bu$ is the vector valued discrete Laplace operator defined by (\ref{lapdis}) for each of its components).
Similarly, we deduce from (\ref{eqenerb}) the discrete  energy balance over the control volume $K$, letting $\theta=1$ in $K$, and $0$ otherwise; this equation reads:
\be\label{nvstksd_nrj}\ba
-\mK \Delta_K T \dsp + \sum_{\edge\in\edgescv, \edge \subset K/L}\Phi_{K,\edge}^{\lambda}(\bu,p)\frac{T_K+T_L} 2
 = \int_K g \d \bx.
\ea\ee
Recall that, for all $K\in\mesh$, and all $\edge\in\edgescv$ such that $\edge\subset\frontomegaun$, the Dirichlet boundary condition (\ref{nvstksd_nrjcld}) is given.
We deduce from (\ref{eqenerb}) the relation imposed by the Neumann boundary condition for the 
thermal flux, letting $\theta_\edge=1$ and $0$ otherwise, for some $\edge\in\edgescv$ with $\edge\subset\frontomegadeux$:
\be\label{nvstksd_nrjcln}\ba
\dsp F_{\cv,\edge}(T) = \int_\edge q_b(\bx) \d s(\bx).
\ea\ee
Note that the above relation is natural, accounting for the fact that $F_{\cv,\edge}(T)$ approximates
the heat flux at the edge $\edge$.
Finally, we write (\ref{eqmasb}) in a given control volume $K$:
\be
\label{nvstksd_mass}
\sum_{\edge\in\edgescv\cap\edgesint} \Phi_{K,\edge}^{\lambda}(\bu,p) = 0.
\ee
{
The stencil of the scheme (\ref{nvstksd_qqtmv}-\ref{nvstksd_mass}) is always determined by that of the diffusion operator.}
\begin{itemize}
	\item[]\em
	As previously mentioned, in the case of orthogonal quadrilateral or parallelepiped grids,
the diffusion flux through an edge $\edge$ gives the classical two-point difference scheme (\ref{diffussion_scheme_superadm}) and the barycenter coordinates are simply calculated by an arithmetic average, $\xedge=(\dcvvedge \xcv+\dcvedge \xcvv)/(\dcvvedge+\dcvedge)$.
Hence, in this case, equations \eqref{nvstksd_qqtmv}-\eqref{nvstksd_mass} then read, for a given $K \in \mesh$:
\begin{align}
\mathbf{Q}_{\Delta u}^{(K)}  + \mathbf{Q}_{\grad p}^{(K)}  +\mathbf{Q}_{u\cdot\grad u}^{(K)}  + \mathbf{Q}_T^{(K)}  & =   \int_K\bff \ \d \bx, \nonumber \\
E_{\Delta T}^{(K)}  + E_{u\cdot\grad T}^{(K)}  & =  \int_K g \d \bx, \nonumber\\
M_{\div u}^{(K)}  & = 0,  \nonumber
\end{align}
with: 
\begin{align}
\mathbf{Q}_{\Delta u}^{(K)} & =  \Pr\left( \sum_{\edge\in\edgescv, \edge \subset K/L} \medge \frac{\bu_\cvv-\bu_\cv} {\dcvedge+\dcvvedge} + \sum_{{\edge\in\edgescv\cap\edgesext}}\medge\frac{\b0-\bu_{K} } {\dcvedge}\right) \nonumber \\ 
\mathbf{Q}_{\grad p}^{(K)}  &= \sum_{\edge\in\edgescv, \edge \subset K/L} 
\medge\frac{\dcvvedge}{\dcvedge+\dcvvedge}(p_\cvv-p_\cv)\ncvedge\nonumber \\ 
\mathbf{Q}_{u\cdot\grad u}^{(K)} &=
 \sum_{\edge\in\edgescv, \edge \subset K/L}\medge
\left(\frac{\dcvvedge\bu_\cv+\dcvedge\bu_\cvv}{\dcvedge+\dcvvedge}\cdot \ncvedge + \lambda_\edge(p_\cv-p_\cvv)\right)\frac{\bu_\cv+\bu_\cvv}{2} \nonumber \\
\mathbf{Q}_T^{(K)}  &=  - \Ra\,\Pr\ \mcv T_K \ez, \nonumber 
\end{align}
\begin{align}
E_{\Delta T}^{(K)}  & =  \sum_{\edge\in\edgescv, \edge \subset K/L} \medge \frac {T_\cvv-T_\cv} {\dcvedge+\dcvvedge} + \sum_{{\edge\in\edgescv\cap\edgesext}}\medge\frac{T_\edge-T_{\cv} } {\dcvedge}\nonumber\\
  E_{u\cdot\grad T}^{(K)} & =  \dsp   \sum_{\edge\in\edgescv, \edge \subset K/L} \medge 
 \left(\frac{\dcvvedge\bu_\cv+\dcvedge\bu_\cvv}{\dcvedge+\dcvvedge} \cdot \ncvedge + \lambda_\edge(p_\cv-p_\cvv)\right)\frac{T_K+T_L} 2,\nonumber
\end{align}
$$ M_{\div u}^{(K)} = \sum_{\edge\in\edgescv, \edge \subset K/ L}\medge \left(  \ncvedge  \frac{\dcvvedge \bu_\cv + \dcvedge \bu_\cvv} {\dcvedge+\dcvvedge} + \lambda_\edge(p_\cv-p_\cvv) \right),$$
and with the Dirichlet boundary (\ref{nvstksd_nrjcld})  applied to $\edge\subset\Gamma_1$ and the Neumann boundary condition on  $\edge\subset\Gamma_2$ being reduced to $ \medge(T_\edge-T_\cv)/\dcvedge= -\int_\edge q_b(\bx) \d s(\bx)$.
\end{itemize}

\subsection{Some mathematical properties}\label{sec-math}
The system of discrete equations (\ref{nvstksd_qqtmv}-\ref{nvstksd_mass})  is a system of non-linear equations. 
The mathematical proof of the existence of at least one solution can be shown in the particular case $T_b = 0$ and $q_b = 0$, which we consider in this section. 
Indeed, in this case, we can show some a priori bounds on $T$ and $\bu$. 
We first let $\theta = T$ in (\ref{eqenerb}). Using the relation
\[
\int_\O \div_\disc^{\lambda}(T,\bu,p) \ P_{\mesh}T \ \d \bx = 0,
\]  
which results from (\ref{eqmasb}), we get
\[
\Vert  \grad_\disc T \Vert_{L^2(\O)^d}^2 = \int_\O g \ P_{\mesh}T \ \d \bx.
\]
Thanks to a discrete Poincar\'e inequality  which follows from \cite[Lemma 5.3]{eym-08-dis}, we get that there exists $C_T$, only depending on the regularity of the mesh and on $g$, but not on the size of the mesh, such that
\[
\Vert  \grad_\disc T \Vert_{L^2(\O)^d} \le C_T.
\]
We then let $\bv = \bu$ in (\ref{eqintmomb}). We get, thanks to (\ref{phi}) and (\ref{eqmasb}),
\[\ba
\dsp \Pr \Vert  \grad_\disc \bu \Vert_{(L^2(\O)^d)^d}^2 +\sum_{\edge\in\edgesint,\edge \subset K/L} \medge \lambda_\edge(p_L - p_K)^2\\
\dsp  = \int_\O (\bff + \Ra\,\Pr \ P_{\mesh}T \ez)\cdot P_{\mesh}\bu \ \d \bx.
\ea\]
Again using the Poincar\'e inequality, we conclude that there  exists $C_\bu$, only depending
on the regularity of the mesh, on $\Ra$, $\Pr$, $\bff$ and $g$, but not on the size of the mesh, such that
\[
\Vert  \grad_\disc \bu \Vert_{(L^2(\O)^d)^d} \le C_\bu.
\]
Hence, using the topological degree method, we can prove the existence of at least one solution.
Moreover, these inequalities are then sufficient to get compactness properties, which show that, from a sequence of discrete solutions with the space step tending to zero, we can extract a converging subsequence,
for suitable norms. Then we can prove that the limit of this subsequence has a sufficient regularity,
in relation with the weak sense provided by (\ref{nstocontfsszero}). 
It is then possible to pass to the
limit on  (\ref{eqenerb}),  (\ref{eqintmomb}) and  (\ref{eqmasb}), using test functions which are 
interpolation of regular ones. We then get that the limit of the converging subsequence satisfies  
(\ref{nstocontfsszero}).

\section{Numerical validation}

Numerical implementation is performed for three dimensional domains. 
We first validate our results on known analytical solutions which allow us to compute the scheme's order of convergence. 
We then turn to a natural convection case which is referenced in the literature. \\
The domains considered are either the unit cube or a circular centered cone. 
In both cases, we use structured (rectangular and non rectangular) meshes, and we denote by $N$ the number of cells in each of the three space directions.\\
The set of non-linear equations (\ref{nvstksd_qqtmv}-\ref{nvstksd_mass}) is solved by an under-relaxed Newton method where the unknowns are the velocity $\bu_\cv$, the pressure $p_\cv$ and the temperature $T_\cv$ and $T_\edge$ for all $\cv \in \mesh$ and $\edge\in\edgesext\cap\frontomegadeux$.
The solutions of the  linear systems are computed with a parallel Generalized Minimal RESidual method provided by the scalable linear solvers package {\em HYPRE} with a preconditioning based on the block Jacobi iLU factorization carried out by the {\em Euclid} library  \cite{HYPRE}.

\subsection{Analytical solutions} \label{sec-analyt}

\begin{figure}[htbp] 
 \centering 
 \begin{minipage}{.5\textwidth}
 \includegraphics[width=\textwidth]{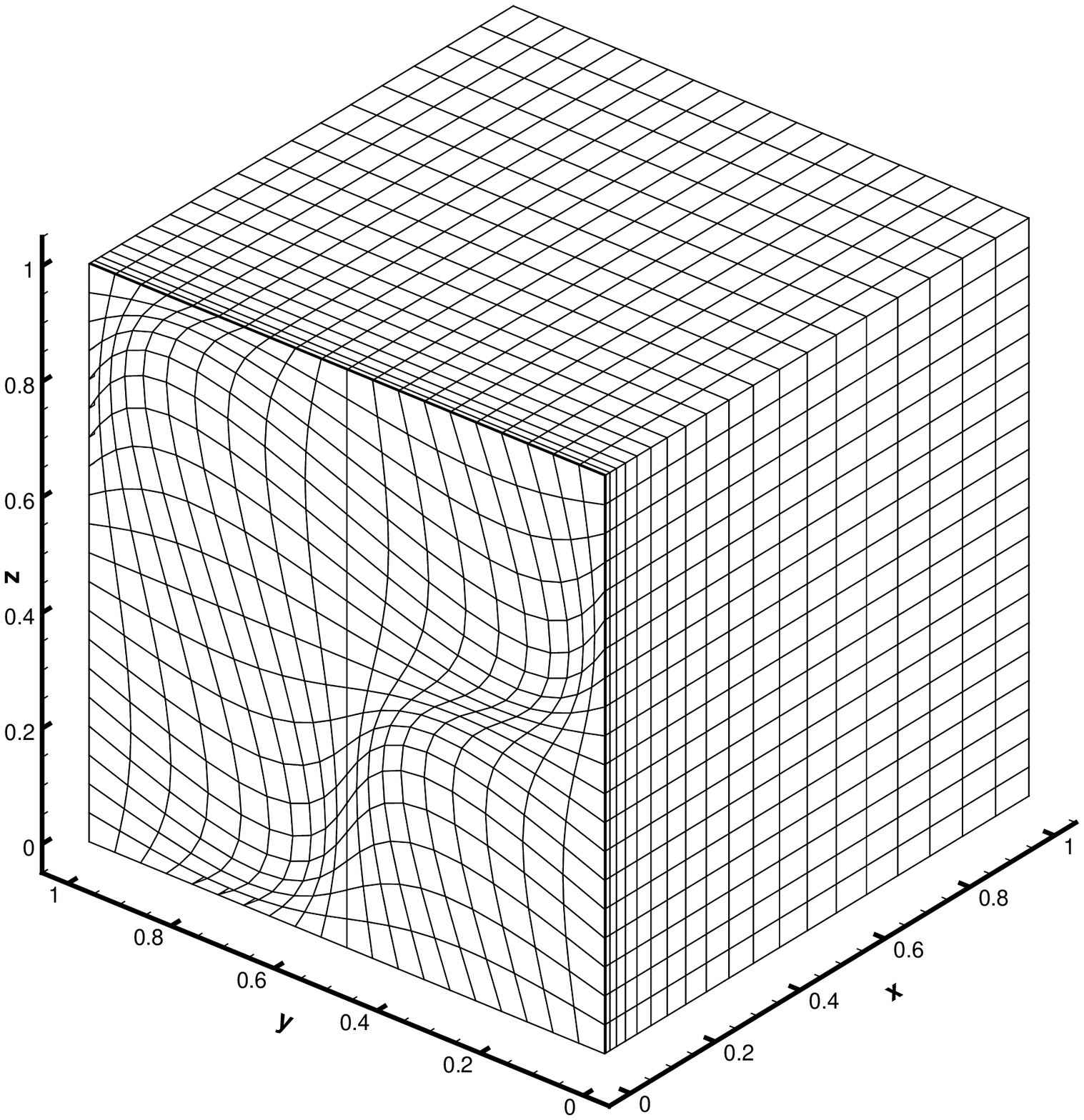}
 \centerline{(a)}
 \end{minipage}\hfill
 \begin{minipage}{.5\textwidth}
 \includegraphics[width=\textwidth]{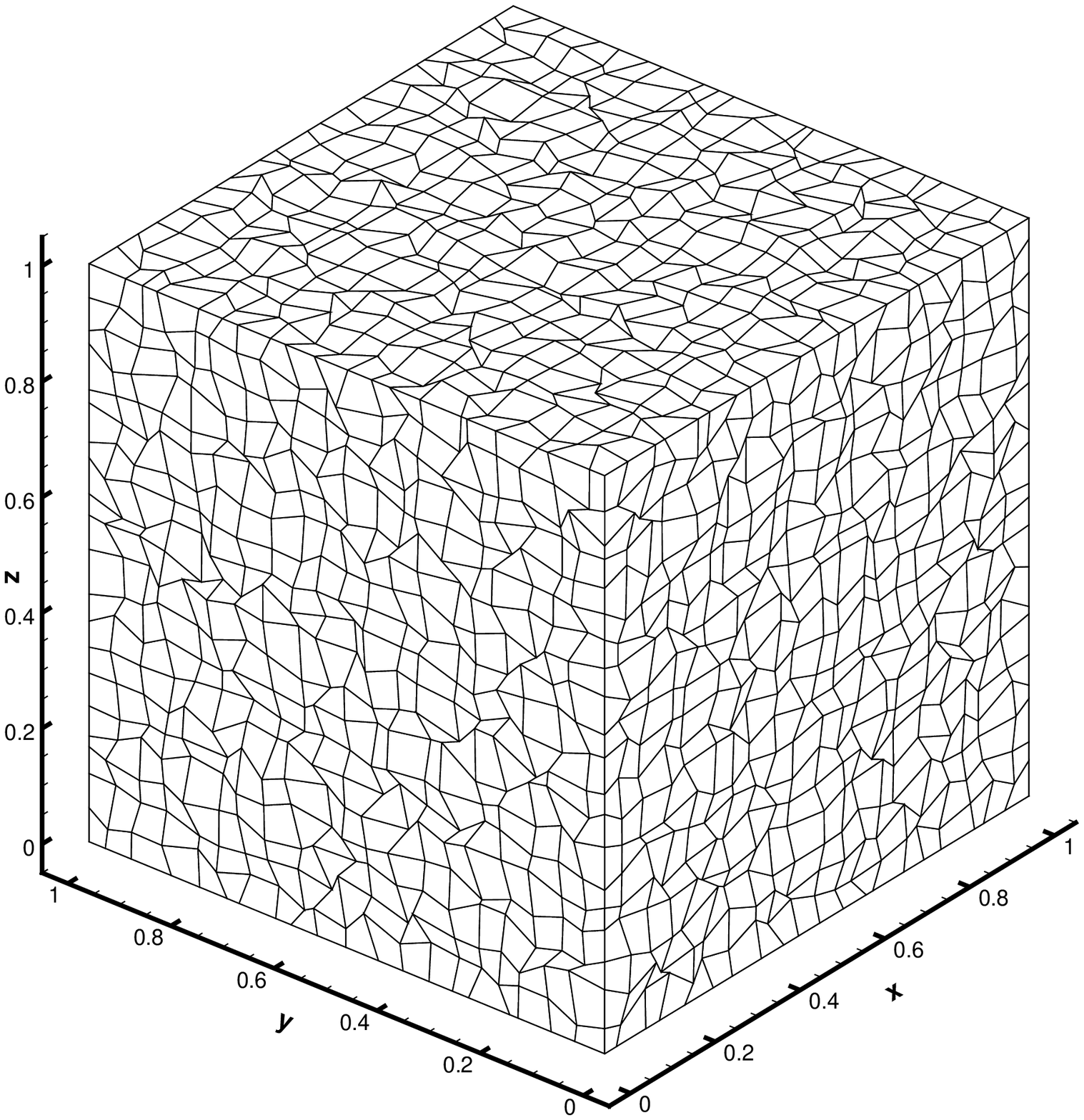}
 \centerline{(b)}
 \end{minipage}\hfill
 \begin{minipage}{.8\textwidth}
 \includegraphics[width=\textwidth]{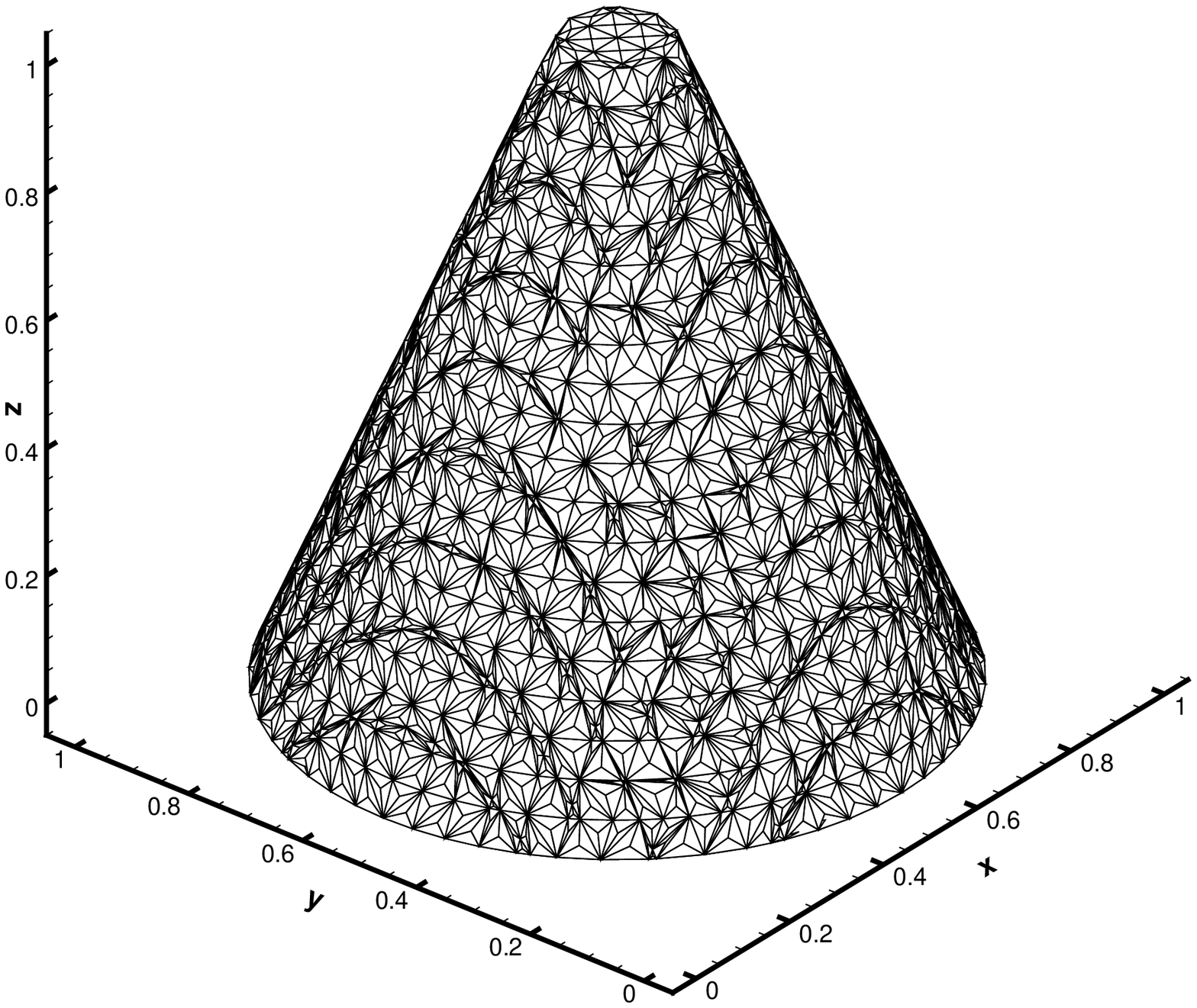}
 \centerline{(c)}
 \end{minipage} 
 \caption[Example of meshes applied for convergence analysis, $N=20$ and $x\equiv x^{(1)}$, $y\equiv x^{(2)}$, $z\equiv x^{(3)}$; (a) Smooth mapping, (b) Random meshes, (c) Truncated cone mesh.]{}
 \label{fig:meshes}
\end{figure}

We consider two closed cavities, cubic or cone-shaped, in which the fluid flow and the heat transfer are known {\em a priori}. 
Let $\pref$, $\buref$ and $\Tref$ be some known pressure, divergence free velocity and  temperature fields; we then compute $\bff$ and $g$ by Eqs.~(\ref{defNVSTKSa},\ref{defNVSTKSc}) where we have set  $\bu \equiv \buref$, $p \equiv\pref $ and $T \equiv\Tref $.\\
For any regular function $\psi$ ($\psi =  (u^{(i)}$, ${i=1,\cdots,3}$, $p$ or $T$), the  scheme's relative accuracy for the usual $L^\infty$, $L^2$ and $H^1$ norms is measured by 
\begin{align}
\epsilon_\infty(\psi) &= \dfrac{\max_{\cv\in\mesh}|\psi_\cv-\psiref(\xcv)|}{\max_{\cv\in\mesh}|\psiref(\xcv)|},\nonumber\\ 
\epsilon_2(\psi)  &= \left(\dfrac{\sum_{\cv\in\mesh}\mcv(\psi_\cv-\psiref(\xcv))^2}{\sum_{\cv\in\mesh}\mcv (\psiref(\xcv))^2}\right)^{1/2},\nonumber\\
\epsilon_{H_1}(\psi)  &= \left(\dfrac{\sum_{\cv\in\mesh}\mcv \vert \grad_\cv \psi -\grad\psiref(\xcv) \vert^2}{\sum_{\cv\in\mesh}\mcv \vert \grad\psiref(\xcv) \vert^2}\right)^{1/2},\nonumber\\
\end{align}
where $\vert \cdot \vert$ denotes the usual Euclidean inner product in $\R^3$.
For each of the above relative error, the scheme's order of convergence is defined by the mean slope of the logarithm of the relative error as a function of the logarithm of the largest cell diameter $\max_{\cv\in\mesh}h_\cv$, the slope being calculated by a least square method.\\
Three different meshes are studied for the unit cubic enclosure.
 Except for the last mesh, the points  $\xcv$ are located at the gravity center of the cells.
\begin{enumerate}
\item The first mesh is an uniform mesh consisting of orthogonal parallelepipeds of size $1/N^3$.
\item The second one (Fig.~\ref{fig:meshes}a)
is constructed by a smooth mapping between the  uniform mesh and the spatial coordinates \cite{Breil2007}.
The vertices $\xs(i,j,k)=\Big( \xsi(l)(i,j,k)\Big)_{l=1,\cdots,3}$ of the elementary distorted cubes are defined by: $\forall (i,j,k)\in N([1,n_1+1])\times N([1,n_2+1])\times N([1,n_3+1])$,
\[
\ba
\xsi(1)(i,j,k)&=&\dsp 1-\cos\left(\frac{\pi(i-1)}{2N}\right)\\
\xsi(2)(i,j,k)&=&\dsp \frac{(j-1)}{N} +0.1\sin\left(\frac{2\pi(j-1)}{N}\right)\sin\left(\frac{2\pi(k-1)}{N}\right)\\
\xsi(3)(i,j,k)&=&\dsp \frac{(k-1)}{N} +0.1\sin\left(\frac{2\pi(j-1)}{N}\right)\sin\left(\frac{2\pi(k-1)}{N}\right).
\ea
\]
\item  The third mesh (Fig.~\ref{fig:meshes}b) is constructed from the first one in the following way: the points $\xcv$ remains at the gravity centers of the basic mesh whereas the vertices $\xs$ of the cells are randomly displaced in each space direction at most of $0.45/N$.
Unlike the previous mesh, which consists in hexahedra with plane faces, the four edges of a face are now not included in a same plane, with the exception of edges which belong to the boundaries of the cubic domain.
Since the consistency of the discrete gradient defined in \eqref{defgradcv} (and therefore of that defined in \eqref{defgrad}) holds if the faces $\sigma$ are plane, we replace each of the non planar faces by two triangular faces.
 \end{enumerate}
The second enclosure, the cone-shaped cavity, is bounded by the lateral surface $((x^{(1)}-0.5)^2+(x^{(2)}-0.5)^2=((6-5x^{(3)})/12 )^2$ for $x^{(3)}\in[0,1]$ and by two plane discs $(x^{(1)}-0.5)^2+(x^{(2)}-0.5)^2\leq 1/4$ for $x^{(3)}=0$ and $(x^{(1)}-0.5)^2+(x^{(2)}-0.5)^2\leq 1/12^2$ for $x^{(3)}=1$.
To make the most of the scheme which locally reduces into a 7-point difference scheme on parallelepiped cells, the mesh consists of cubes which were cut to match the lateral curved boundary.
Thus, the mesh error tends quadratically to zero  with respect to the mesh size.
To avoid too large differences of volume sizes between adjacent cells that may deteriorate the numerical accuracy, the boundary cells having a volume less than $0.1/(n_i)^d$ are merged into adjacent cells.

\begin{figure}[htbp] 
 \centering 
 \begin{minipage}{.5\textwidth}
 \includegraphics[height=\textwidth,angle=-90]{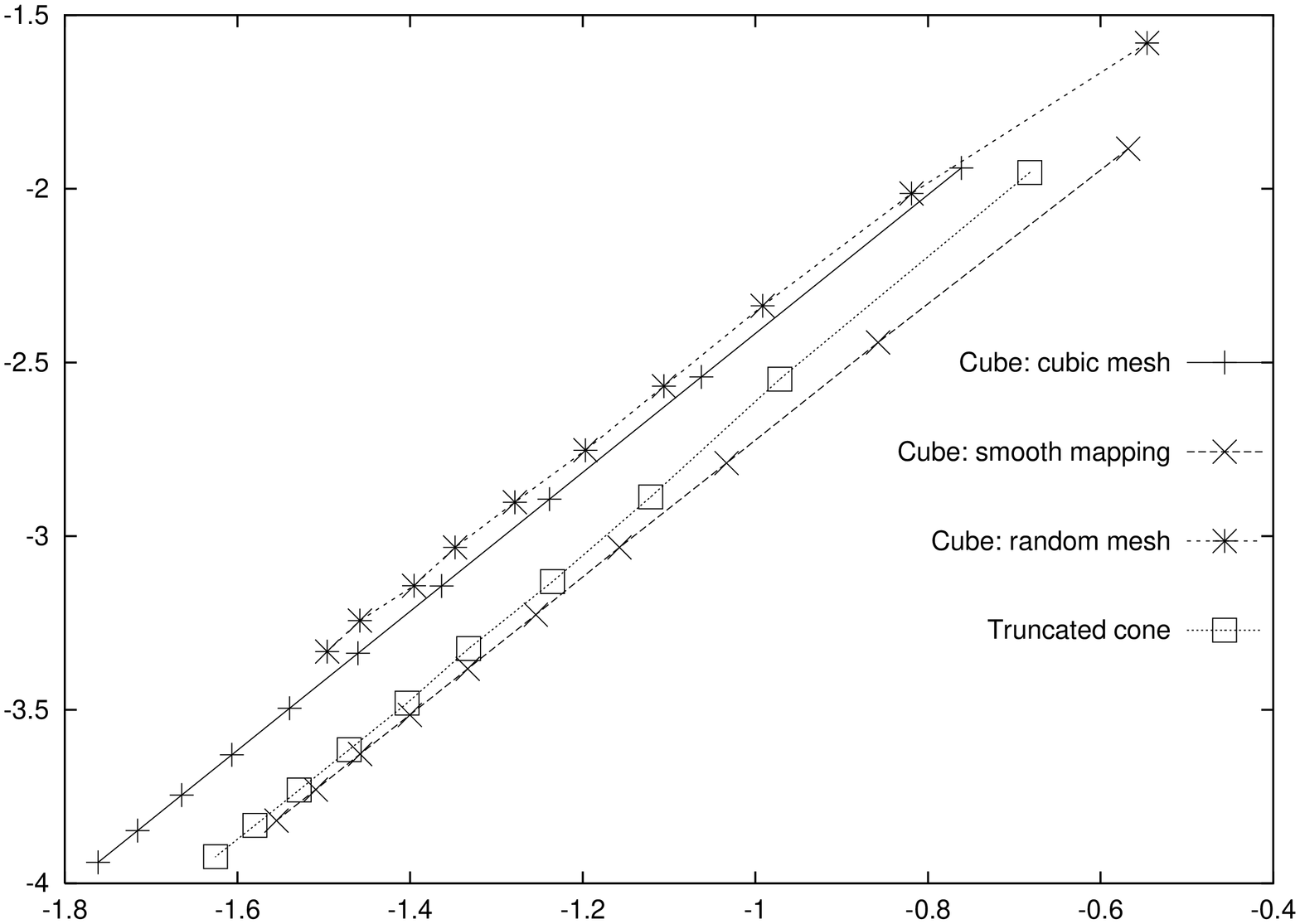}
 \centerline{(a)}
 \end{minipage}\hfill
 \begin{minipage}{.5\textwidth}
 \includegraphics[height=\textwidth,angle=-90]{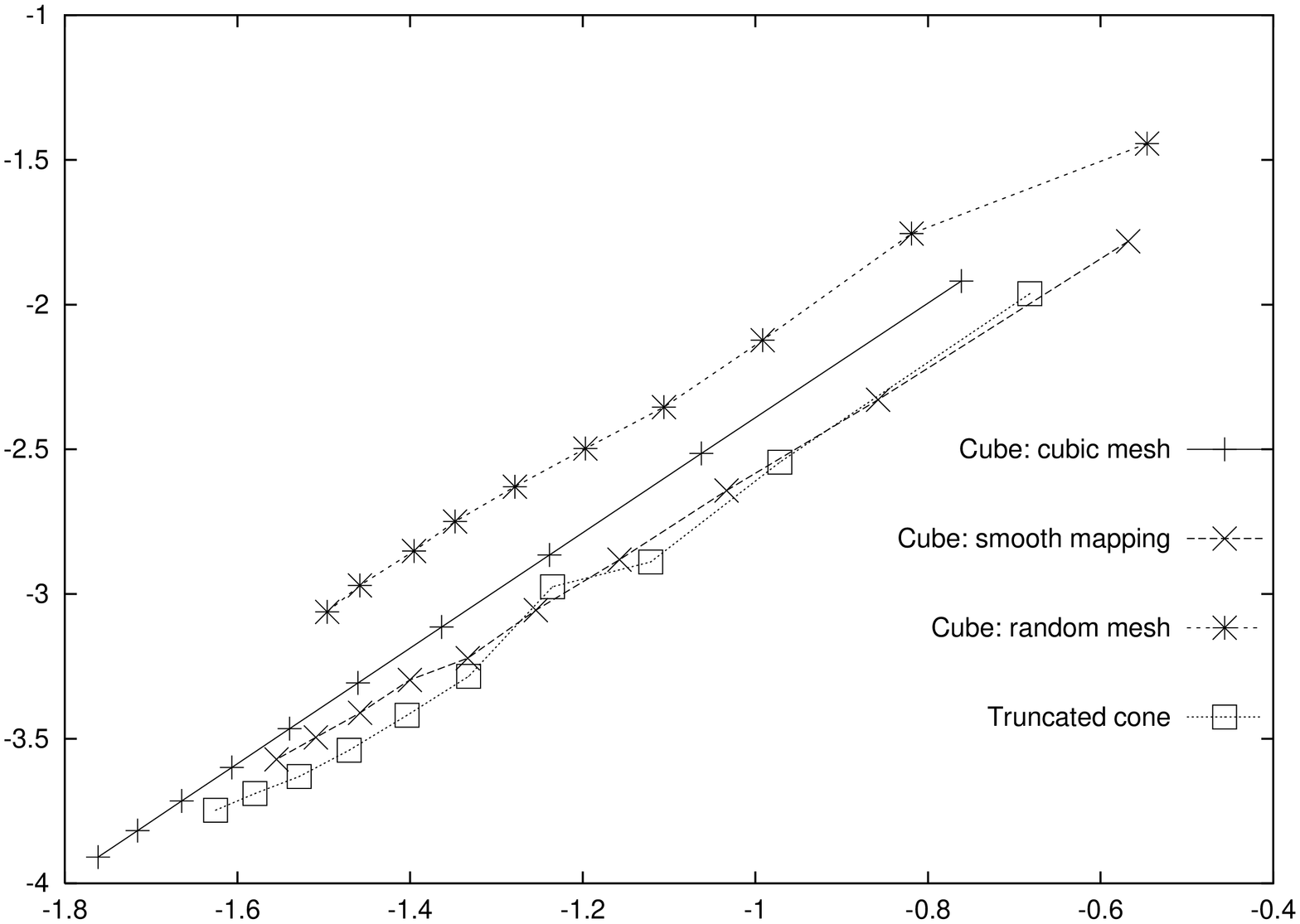}
 \centerline{(b)}
 \end{minipage}\hfill
 \begin{minipage}{.5\textwidth}
 \includegraphics[height=\textwidth,angle=-90]{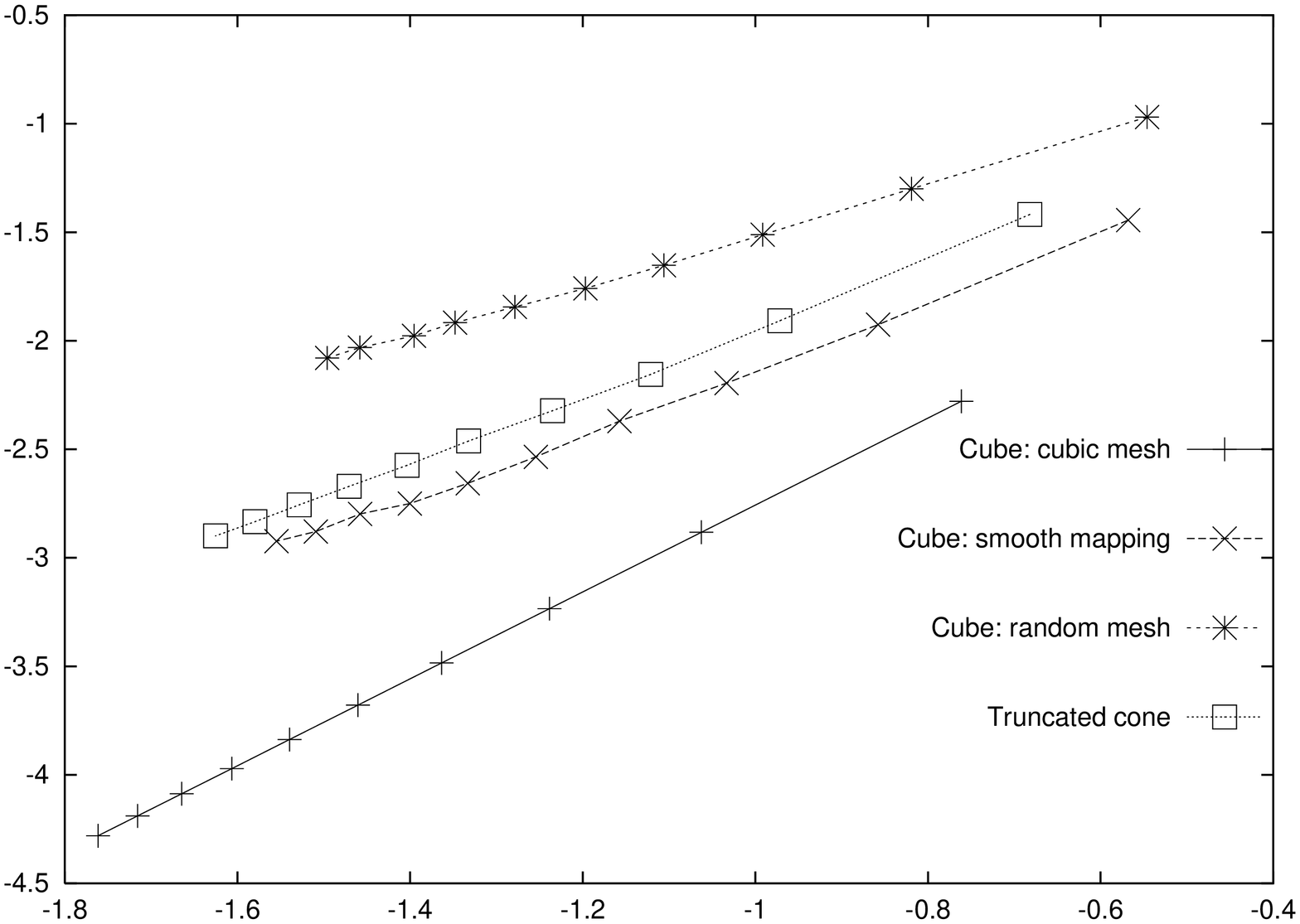}
 \centerline{(c)}
 \end{minipage}
 \caption[
 Solution accuracy of the Poisson problem as a function of $\log(\max_{\cv\in\mesh}h_\cv)$, for $10\leq N\leq 100$ and different geometries and meshes; (a) $\log(\epsilon_2(T))$, (b) $\log(\epsilon_\infty(T))$, (c) $\log(\epsilon_{H_1}(T))$.]{}
 \label{fig:accuracy_diff}
\end{figure}

We are first interested in the Poisson problem for a scalar variable (Eq. (\ref{nvstksd_nrj}) with $\bu = \b0$) where $\frontomegaun=\frontomega$ and $\frontomegadeux=\emptyset$.
It was first checked that the errors obtained with a linear analytical solution on the different meshes and cavities are of the order of the computer accuracy, even for the coarsest grids.
The next analytical test consists in choosing the reference solution \(\Tref(x^{(1)},x^{(2)},x^{(3)})=\sin(\pi x^{(1)})\cos(\pi x^{(2)})\cos(\pi x^{(3)})\) with appropriate Dirichlet boundary conditions (Fig.~\ref{fig:accuracy_diff}a-c).
The orders of convergence are reported in table~(\ref{tab:slopes_diff}).
The accuracy of the scheme is close to $2$ when considering the $L^2$-norm and it slightly decreases with the $L^\infty$-norm but always remains larger than $1.50$.
The order of convergence for the gradients ($H^1$-norm) is larger than $1$. \\
Next, we examine the convergence behavior of the isothermal Navier-Stokes equations by setting 
$\buref(\bx)=\grad \wedge \sum_{i=1}^d (4x^{(1)}(x^{(1)}-1))^3 (4x^{(2)}(x^{(2)}-1))^4\break (4x^{(3)}(x^{(3)}-1))^5\ei$ and $\pref(\bx)= \cos(\pi x^{(1)})\cos(\pi x^{(2)})\cos(\pi x^{(3)})$ 
in (Eq.~\ref{defNVSTKSa}) with $\Pr=1$ and $\Ra=0$ 
(note that the dimensionless writing of the equations is meaningless because the current reference velocity is related to the thermal diffusivity which has no reason to appear for isothermal problems. 
Another velocity reference should be used, based on the viscous diffusivity).
Table~(\ref{tab:slopes_nvstks}) 
indicates that the convergence rates of the velocity components are larger than $1.90$ on the three finer meshes when the relative error is based on the $L^2$-norm and first order accurate for the pressure for distorted meshes.
In accordance with the diffusion problem when the $L^\infty$-norm is used, the orders of convergence slightly decrease for the velocity but a convergence rate larger than $1.60$ is still measured.
The convergence rates of the gradients are better than the expected first order.
Unsurprisingly, the $L^\infty$ and $H^1$-norms of the pressure do not tend to zero with the mesh size because it simply appears in the momentum equation as Lagrangian multiplier of the mass equation.
Thus the only guaranteed convergence for the pressure is based on the $L^2$-norm.

\subsection{Natural convection problem\label{Natural_convection_problem}}

We consider an air filled unit-cubic enclosure with isolated walls except the two face to face vertical isothermal surfaces at $x^{(1)}=0$ and $1$. 
The governing fluid flow equations are solution of system (\ref{defNVSTKS}) and (\ref{diff_cls}) with $\bff(\bx)=\b0$, $g(\bx)=0$, $q_b(\bx)=0$, $T(0,x^{(2)},x^{(3)})=-0.5$ and $T(1,x^{(2)},x^{(3)})=0.5$.
The Prandtl and Rayleigh numbers are fixed to $\Pr=0.71$ and $\Ra=10^7$ and the stabilization parameter is chosen equal to $\lambda=10^{-8}$ in the mass equation.
 Because very small boundary layers take place along the walls, the vertices are located at the Gauss-Lobatto points and the collocation points $\xcv$ at the gravity centers of the cells.
To also study the effect of non-cubic meshes, the coordinates of the previous defined vertices are randomly moved, for each direction, of a magnitude at most equal to $0.45$ the size of the cell in this direction.

Table~(\ref{cvn}) 
presents the maxima of the velocity components, the average Nusselt number on the isothermal walls and their relative differences with respect to reference data \cite{Tric00}.
For the cubic meshes, the relative differences with respect to reference data seem to tend to zero:
for the finer mesh, our results depart from less than $1\%$ with respect to the reference values.
Although the solutions are less accurate for the shaken meshes, a convergent behavior is also observed.
The relatively large differences measured on $u_\infty^{(2)}$, in comparison with the other velocity components, are probably the result of the flow shape in which the main motion occurs in the $(\ex,\ez)$-planes with a small secondary flow in the transverse planes.
It is also interesting to note the good accuracy of the average Nusselt numbers.
This accuracy is essentially obtained thanks to the use of the stabilized mass flux $\Phi_{K,\edge}^{\lambda}(\bu,p)$ in the heat transport expression which
ensures the conservation of the average heat flux balance on the boundaries of the cavity.


\section{Conclusion}

In this paper  we presented a new scheme which is well suited for the simulation of incompressible viscous flows on irregular and non-conforming grids.
This possibility seems to open a large field of new applications (grid refinement as a function
of an {\em a posteriori} error computation, free boundaries, \ldots).
We emphasize that the convergence of the scheme may be proven mathematically, and that the obtained numerical results are accurate. 
Although we presented this scheme in the steady case, its extension to transient regimes is straightforward. In this latter case, one should consider optimizing the  linear solving step by  using suitable projection algorithms.




\bibliographystyle{elsart-num-sort}

\newpage

\listoftables

\newpage

\listoffigures

\newpage
\begin{table}[htbp]
 \centering 
\[
\begin{array}{||c|c|c|c|c||}
\hline\hline
 &\multicolumn{3}{|c|}{\hbox{Cubic enclosure}}&\hbox{Truncated conic}\\
 &\hbox{cubic meshes}&\hbox{smooth meshes}&\hbox{random meshes}&\hbox{enclosure}\\
 \hline
L^2\hbox{-norm}& 2\ (2) &1.96\ (1.97) &1.87\ (1.93) & 2.09 \ (2.03)\\
L^\infty\hbox{-norm}& 1.99\ (2) &1.81\ (1.68) &1.74\ (1.88) & 1.92 \ (1.58)\\
H^1\hbox{-norm} & 2\ (2) &1.50\ (1.28) &1.16\ (1.07) & 1.55 \ (1.46)\\
\hline\hline
\ea
\]
\caption[Orders of convergence for the Poisson problem, $10\leq N\leq 100$ ($50\leq N\leq 100$).]{\label{tab:slopes_diff}}
\end{table}

\newpage

\begin{table}[htbp]
 \centering 
\[
\begin{array}{||c|c|c|c|c||}
\hline\hline
 & & \hbox{Cubic meshes} &\hbox{Smooth meshes} &\hbox{Random meshes} \\
 \hline
L^2 &\ui(1) & 2\ (2) & 1.94\ (1.99) & 1.82\ (1.95)\\ 
 &\ui(2) & 1.99\ (2) & 1.96\ (1.99)& 1.84\ (1.94)\\ 
 &\ui(3) & 2\ (2) & 1.95\ (1.98)& 1.87\ (1.94)\\ 
 &p & 2\ (2) & 1.09\ (0.76)& 0.85\ (0.93)\\ 
\hline
L^\infty &\ui(1) & 2\ (2) & 1.67\ (2.02) & 1.61\ (1.84)\\ 
 &\ui(2) & 1.75\ (1.78)& 1.39\ (1.64)& 1.66\ (1.72)\\ 
 &\ui(3) & 1.88\ (1.74) & 1.55\ (1.59)& 1.73\ (1.86)\\ 
 &p & 1.72\ (1.93) & -- & --\\ 
\hline
H^1 &\ui(1) & 1.98\ (2) & 1.80\ (1.85) & 1.43\ (1.31)\\ 
 &\ui(2) & 2\ (2) & 1.81\ (1.84)& 1.41\ (1.27)\\ 
 &\ui(3) & 1.95\ (1.99) & 1.76\ (1.83)& 1.40\ (1.25)\\ 
 &p & 1.91\ (1.97) & -- & --\\ 
\hline\hline
\ea
\]
\caption[Orders of convergence for the isothermal Navier-Stokes problem, $10\leq N \leq 60$ ($40\leq N\leq 60$).]{\label{tab:slopes_nvstks}}
\end{table}

\newpage

\begin{table}[htbp]\label{cvn}
 \centering 
\[
\begin{array}{||c|c|c|c|c|c|c|c|c|c||}
\hline\hline
 \hbox{Mesh types} & N &u^{(1)}_\infty& e(u^{(1)}_\infty)&u^{(2)}_\infty& e(u^{(2)}_\infty)&u^{(3)}_\infty& e(u^{(3)}_\infty)&\Nu&e(\Nu)\\
 \hline
 &20&333,23&-13\% &70,959&-15\% &767,01&-0.15\%&16,380&0.23\%\\
 &30&371,89&-3.1\%&79,105&-5.1\%&761,11&-0.91\%&16,366&0.14\%\\
\hbox{cubic}
 &40&377,71&-1.6\% &81,097&-2.7\%&761,15&-0.91\%&16,361&0.11\% \\
 &50&380,19&-0.95\%&82,234&-1.4\% &767,25&-0.11\% &16,357&0.086\%\\
 &60&380,47&-0.88\%&82,615&-0.93\%&767,90&-0.031\%&16,353&0.065\%\\
 \hline
 &20&497,85&30\%&363,48&340\%&869,95&13\%&16,023&-2.0\%\\
\hbox{Random}
 &30&419,43&9.3\%&276,63&230\%&777,25&1.2\%&16,198&-0.88\%\\
 &40&400,18&4.3\%&151,27&81\% &779,93&1.5\%&16,259&-0.51\%\\
 \hline
 \hbox{\cite{Tric00}}& & 383,8357&0\%& 83,3885&0\%&768,1393&0\%&16,3427&0\%\\
 \hline\hline
\ea
\]
\caption[Maxima of the velocity components $u^{(i)}_\infty = \max_{\cv\in \mesh} |u^{(i)}|$ ($i= 1,\cdots,3$), average Nusselt number $\Nu=\int_0^1\int_0^1(\grad T\cdot {\bf n})_{x_1=0}\,\d x_2\d x_3$ and their relative differences to reference values \cite{Tric00}, $e(\psi)=\psi_\infty/\psi_{\rm ref}-1$, $\psi\in\{u^{(1)}_\infty,u^{(2)}_\infty,u^{(3)}_\infty,\Nu\}$ for the natural convection problem with $\Ra=10^7$ and $\Pr=0.71$.]{}
\end{table}

\newpage

\begin{figure}[htb]
\begin{center}
\resizebox{\textwidth}{!}{\input{mesh_new.pstex_t}}
\caption[Cell $\cv$ and neighbors for a mesh of a two-dimensional ($d=2$) domain $\Omega$: $\edge$ an interior edge ($ \edge \subset K/L$), $\edgep$ a boundary edge ($ \edgep \in \edges_\cv$), 
$\mcv$ the measure of the cell $\cv$,
$\xcv$ the collocation point, 
$\xedge$ the barycenter of $\edge$ of measure $\medge$, 
$\dcvedge$ the Euclidean distance between $\xcv$ and $\edge$, 
$\ncvedge$ the unit outward vector normal to $\edge$ 
and $C_{K,\edge}$ the cone with vertex $\xcv$ and basis $\edge$ (similar notations apply for the edge $\edgep$).]{}
\label{figmesh}
\end{center}
\end{figure}
\newpage

\begin{figure}[htb]
\begin{center}
\resizebox{\textwidth}{!}{\input{mesh_cart.pstex_t}}
\caption[Cell $\cv$ and  neighbors for a mesh of a two-dimensional ($d=2$) domain $\Omega$ meshed with an uniform rectangular grid with step size $h_x$ in the $x$ direction and $h_y$ in the $y$ direction: $\edge=\edge_{i,j-1/2}$ and $\edgep=\edge_{i+1/2,j}$ two interior edges with $\edge_{i,j-1/2} \subset K_{i,j}/K_{i,j-1}$ and $ \edge_{i+1/2,j}\subset K_{i,j}/K_{i+1,j}$, 
$\mcv = h_x h_y$ the measure of any cell $\cv$,
$\xcv$ the collocation point (center of gravity), 
$\xedge$ the barycenter of $\edge$ of measure $\medge$, 
$\dcvedge= {h_x}/ 2$ the Euclidean distance between $\xcv$ and $\edge$,
$\ncvedge$ the unit outward vector normal to $\edge$ 
and $C_{K,\edge}$ the cone with vertex $\xcv$ and basis $\edge$.]{}
\label{figmeshcart}
\end{center}
\end{figure}

\newpage

\begin{figure}[htb]
\begin{center}
\begin{minipage}{.48\textwidth}
\resizebox{\textwidth}{!}{\input{gradcart.pstex_t}}
\centerline{(a)}
\end{minipage}
\begin{minipage}{.48\textwidth}
\resizebox{\textwidth}{!}{\input{gradcartmodif.pstex_t}}
\centerline{(b)}
\end{minipage}
\caption[Gradient of the element $1_{\cv_{i,j}}$ of $\dset$ for an uniform rectangular grid with step size $h_x$ in the $x$ direction and $h_y$ in the $y$ direction; (a) centered gradient, (b) stabilized gradient.]{}
\label{figgradcart}
\end{center}
\end{figure}

\newpage

\newpage

\newpage

\end{document}

%% file: mesh_new.pstex_t
\begin{picture}(0,0)%
\includegraphics{mesh_new.pstex}%
\end{picture}%
\setlength{\unitlength}{3947sp}%
\begingroup\makeatletter\ifx\SetFigFont\undefined%
\gdef\SetFigFont#1#2#3#4#5{%
  \reset@font\fontsize{#1}{#2pt}%
  \fontfamily{#3}\fontseries{#4}\fontshape{#5}%
  \selectfont}%
\fi\endgroup%
\begin{picture}(5987,5396)(2078,-7084)
\put(5752,-4586){\makebox(0,0)[lb]{\smash{{\SetFigFont{12}{14.4}{\familydefault}{\mddefault}{\updefault}{\color[rgb]{0,0,0}$\xedge$}%
}}}}
\put(5702,-3859){\rotatebox{52.0}{\makebox(0,0)[lb]{\smash{{\SetFigFont{12}{14.4}{\familydefault}{\mddefault}{\updefault}{\color[rgb]{0,0,0}$d_{K,\edge}$}%
}}}}}
\put(4299,-3458){\rotatebox{9.0}{\makebox(0,0)[b]{\smash{{\SetFigFont{12}{14.4}{\familydefault}{\mddefault}{\updefault}{\color[rgb]{0,0,0}$h_{K}$}%
}}}}}
\put(5266,-5306){\rotatebox{55.0}{\makebox(0,0)[lb]{\smash{{\SetFigFont{12}{14.4}{\familydefault}{\mddefault}{\updefault}{\color[rgb]{0,0,0}$\edge$}%
}}}}}
\put(5489,-1835){\rotatebox{340.0}{\makebox(0,0)[lb]{\smash{{\SetFigFont{12}{14.4}{\familydefault}{\mddefault}{\updefault}{\color[rgb]{0,0,0}$\medgep$}%
}}}}}
\put(5295,-2273){\rotatebox{341.0}{\makebox(0,0)[lb]{\smash{{\SetFigFont{12}{14.4}{\familydefault}{\mddefault}{\updefault}{\color[rgb]{0,0,0}$\edgep$}%
}}}}}
\put(5947,-2500){\makebox(0,0)[lb]{\smash{{\SetFigFont{12}{14.4}{\familydefault}{\mddefault}{\updefault}{\color[rgb]{0,0,0}$\xedgep$}%
}}}}
\put(5114,-3610){\makebox(0,0)[rb]{\smash{{\SetFigFont{12}{14.4}{\familydefault}{\mddefault}{\updefault}{\color[rgb]{0,0,0}$\xcv$}%
}}}}
\put(5735,-5541){\rotatebox{55.0}{\makebox(0,0)[lb]{\smash{{\SetFigFont{12}{14.4}{\familydefault}{\mddefault}{\updefault}{\color[rgb]{0,0,0}$\medge$}%
}}}}}
\put(6831,-5322){\makebox(0,0)[b]{\smash{{\SetFigFont{10}{12.0}{\familydefault}{\mddefault}{\updefault}{\color[rgb]{0,0,0}$\cvv$}%
}}}}
\put(7121,-3900){\makebox(0,0)[b]{\smash{{\SetFigFont{10}{12.0}{\familydefault}{\mddefault}{\updefault}{\color[rgb]{0,0,0}$\cvvvv$}%
}}}}
\put(5529,-2937){\rotatebox{332.0}{\makebox(0,0)[lb]{\smash{{\SetFigFont{12}{14.4}{\familydefault}{\mddefault}{\updefault}{\color[rgb]{0,0,0}$d_{K,\edgep}$}%
}}}}}
\put(6736,-4443){\makebox(0,0)[lb]{\smash{{\SetFigFont{12}{14.4}{\familydefault}{\mddefault}{\updefault}{\color[rgb]{0,0,0}$\ncvedge$}%
}}}}
\put(6701,-1984){\makebox(0,0)[lb]{\smash{{\SetFigFont{12}{14.4}{\familydefault}{\mddefault}{\updefault}{\color[rgb]{0,0,0}$\ncvedgep$}%
}}}}
\put(7127,-6129){\makebox(0,0)[lb]{\smash{{\SetFigFont{14}{16.8}{\familydefault}{\mddefault}{\updefault}{\color[rgb]{0,0,0}$\Gamma$}%
}}}}
\put(3676,-5386){\makebox(0,0)[b]{\smash{{\SetFigFont{10}{12.0}{\familydefault}{\mddefault}{\updefault}{\color[rgb]{0,0,0}$\cvvv$}%
}}}}
\put(4801,-6211){\makebox(0,0)[b]{\smash{{\SetFigFont{14}{16.8}{\familydefault}{\mddefault}{\updefault}{\color[rgb]{0,0,0}$\Omega$}%
}}}}
\put(4051,-2836){\makebox(0,0)[b]{\smash{{\SetFigFont{12}{14.4}{\familydefault}{\mddefault}{\updefault}{\color[rgb]{0,0,0}$\cv$}%
}}}}
\put(5504,-4253){\makebox(0,0)[b]{\smash{{\SetFigFont{12}{14.4}{\familydefault}{\mddefault}{\updefault}{\color[rgb]{0,0,0}$C_{K,\edge}$}%
}}}}
\put(5330,-2745){\makebox(0,0)[b]{\smash{{\SetFigFont{12}{14.4}{\familydefault}{\mddefault}{\updefault}{\color[rgb]{0,0,0}$C_{K,\edgep}$}%
}}}}
\end{picture}%

%% file: mesh_cart.pstex_t
\begin{picture}(0,0)%
\includegraphics{mesh_cart.pstex}%
\end{picture}%
\setlength{\unitlength}{3947sp}%
\begingroup\makeatletter\ifx\SetFigFont\undefined%
\gdef\SetFigFont#1#2#3#4#5{%
  \reset@font\fontsize{#1}{#2pt}%
  \fontfamily{#3}\fontseries{#4}\fontshape{#5}%
  \selectfont}%
\fi\endgroup%
\begin{picture}(8669,5594)(1779,-6758)
\put(9751,-3361){\makebox(0,0)[lb]{\smash{{\SetFigFont{12}{14.4}{\familydefault}{\mddefault}{\updefault}{\color[rgb]{0,0,0}$\ncvedgep$}%
}}}}
\put(7801,-1861){\makebox(0,0)[lb]{\smash{{\SetFigFont{12}{14.4}{\familydefault}{\mddefault}{\updefault}{\color[rgb]{0,0,0}$d_{K,\edgep}=\frac {h_x} 2$}%
}}}}
\put(6901,-4711){\makebox(0,0)[lb]{\smash{{\SetFigFont{12}{14.4}{\familydefault}{\mddefault}{\updefault}{\color[rgb]{0,0,0}$\ncvedge$}%
}}}}
\put(6451,-1936){\makebox(0,0)[b]{\smash{{\SetFigFont{12}{14.4}{\familydefault}{\mddefault}{\updefault}{\color[rgb]{0,0,0}$\cv= K_{i,j}$}%
}}}}
\put(7501,-5011){\makebox(0,0)[lb]{\smash{{\SetFigFont{12}{14.4}{\familydefault}{\mddefault}{\updefault}{\color[rgb]{0,0,0}$\medge = {h_x} $}%
}}}}
\put(6826,-4186){\makebox(0,0)[lb]{\smash{{\SetFigFont{12}{14.4}{\familydefault}{\mddefault}{\updefault}{\color[rgb]{0,0,0}$\edge= \edge_{i,j-1/2}$}%
}}}}
\put(4801,-2236){\makebox(0,0)[b]{\smash{{\SetFigFont{10}{12.0}{\familydefault}{\mddefault}{\updefault}{\color[rgb]{0,0,0}$\cvvv   =K_{i-1,j}$}%
}}}}
\put(3826,-3061){\rotatebox{90.0}{\makebox(0,0)[lb]{\smash{{\SetFigFont{12}{14.4}{\familydefault}{\mddefault}{\updefault}{\color[rgb]{0,0,0}$ {h_y}$}%
}}}}}
\put(7801,-5761){\makebox(0,0)[b]{\smash{{\SetFigFont{10}{12.0}{\familydefault}{\mddefault}{\updefault}{\color[rgb]{0,0,0}$\cvv=K_{i,j-1}$}%
}}}}
\put(10276,-2011){\rotatebox{270.0}{\makebox(0,0)[lb]{\smash{{\SetFigFont{12}{14.4}{\familydefault}{\mddefault}{\updefault}{\color[rgb]{0,0,0}$\medgep=h_y$}%
}}}}}
\put(7876,-3886){\makebox(0,0)[lb]{\smash{{\SetFigFont{12}{14.4}{\familydefault}{\mddefault}{\updefault}{\color[rgb]{0,0,0}$\xedge= {\bf x}_{i,j-1/2}$}%
}}}}
\put(9826,-1561){\rotatebox{270.0}{\makebox(0,0)[lb]{\smash{{\SetFigFont{12}{14.4}{\familydefault}{\mddefault}{\updefault}{\color[rgb]{0,0,0}$\edgep= \edge_{i+1/2,j}$}%
}}}}}
\put(6270,-2576){\rotatebox{270.0}{\makebox(0,0)[lb]{\smash{{\SetFigFont{12}{14.4}{\familydefault}{\mddefault}{\updefault}{\color[rgb]{0,0,0}$d_{K,\edge}=\frac {h_y} 2$}%
}}}}}
\put(7856,-2603){\makebox(0,0)[lb]{\smash{{\SetFigFont{12}{14.4}{\familydefault}{\mddefault}{\updefault}{\color[rgb]{0,0,0}$\xcv= {\bf x}_{i,j}$}%
}}}}
\put(9519,-2938){\makebox(0,0)[rb]{\smash{{\SetFigFont{12}{14.4}{\familydefault}{\mddefault}{\updefault}{\color[rgb]{0,0,0}$\xedgep={\bf x}_{i+1/2,j}$}%
}}}}
\put(7801,-6586){\makebox(0,0)[b]{\smash{{\SetFigFont{10}{12.0}{\familydefault}{\mddefault}{\updefault}{\color[rgb]{0,0,0}$i$}%
}}}}
\put(3976,-6586){\makebox(0,0)[b]{\smash{{\SetFigFont{10}{12.0}{\familydefault}{\mddefault}{\updefault}{\color[rgb]{0,0,0}$i-1$}%
}}}}
\put(2176,-5236){\makebox(0,0)[b]{\smash{{\SetFigFont{10}{12.0}{\familydefault}{\mddefault}{\updefault}{\color[rgb]{0,0,0}$j-1$}%
}}}}
\put(2176,-2761){\makebox(0,0)[b]{\smash{{\SetFigFont{10}{12.0}{\familydefault}{\mddefault}{\updefault}{\color[rgb]{0,0,0}$j$}%
}}}}
\put(7927,-3480){\makebox(0,0)[b]{\smash{{\SetFigFont{12}{14.4}{\familydefault}{\mddefault}{\updefault}{\color[rgb]{0,0,0}$C_{K,\edge}$}%
}}}}
\put(9063,-2312){\makebox(0,0)[b]{\smash{{\SetFigFont{12}{14.4}{\familydefault}{\mddefault}{\updefault}{\color[rgb]{0,0,0}$C_{K,\edgep}$}%
}}}}
\put(7111,-3153){\rotatebox{33.0}{\makebox(0,0)[b]{\smash{{\SetFigFont{12}{14.4}{\familydefault}{\mddefault}{\updefault}{\color[rgb]{0,0,0}$h_{K}=\sqrt{h_x^2+h_y^2}$}%
}}}}}
\end{picture}%

%% file: gradcart.pstex_t
\begin{picture}(0,0)%
\includegraphics{gradcart.pstex}%
\end{picture}%
\setlength{\unitlength}{3947sp}%
\begingroup\makeatletter\ifx\SetFigFont\undefined%
\gdef\SetFigFont#1#2#3#4#5{%
  \reset@font\fontsize{#1}{#2pt}%
  \fontfamily{#3}\fontseries{#4}\fontshape{#5}%
  \selectfont}%
\fi\endgroup%
\begin{picture}(3120,4014)(6043,-4963)
\put(6751,-4936){\makebox(0,0)[b]{\smash{{\SetFigFont{6}{7.2}{\familydefault}{\mddefault}{\updefault}{\color[rgb]{0,0,0}$i-1$}%
}}}}
\put(7651,-4936){\makebox(0,0)[b]{\smash{{\SetFigFont{6}{7.2}{\familydefault}{\mddefault}{\updefault}{\color[rgb]{0,0,0}$i$}%
}}}}
\put(8551,-4936){\makebox(0,0)[b]{\smash{{\SetFigFont{6}{7.2}{\familydefault}{\mddefault}{\updefault}{\color[rgb]{0,0,0}$i+1$}%
}}}}
\put(6151,-1711){\makebox(0,0)[b]{\smash{{\SetFigFont{6}{7.2}{\familydefault}{\mddefault}{\updefault}{\color[rgb]{0,0,0}$j+1$}%
}}}}
\put(6151,-4111){\makebox(0,0)[b]{\smash{{\SetFigFont{6}{7.2}{\familydefault}{\mddefault}{\updefault}{\color[rgb]{0,0,0}$j-1$}%
}}}}
\put(6151,-2911){\makebox(0,0)[b]{\smash{{\SetFigFont{6}{7.2}{\familydefault}{\mddefault}{\updefault}{\color[rgb]{0,0,0}$j$}%
}}}}
\put(6751,-1711){\makebox(0,0)[b]{\smash{{\SetFigFont{8}{9.6}{\familydefault}{\mddefault}{\updefault}{\color[rgb]{0,0,0}${\bf0}$}%
}}}}
\put(8551,-1711){\makebox(0,0)[b]{\smash{{\SetFigFont{8}{9.6}{\familydefault}{\mddefault}{\updefault}{\color[rgb]{0,0,0}${\bf0}$}%
}}}}
\put(8551,-4111){\makebox(0,0)[b]{\smash{{\SetFigFont{8}{9.6}{\familydefault}{\mddefault}{\updefault}{\color[rgb]{0,0,0}${\bf0}$}%
}}}}
\put(6751,-4111){\makebox(0,0)[b]{\smash{{\SetFigFont{8}{9.6}{\familydefault}{\mddefault}{\updefault}{\color[rgb]{0,0,0}${\bf0}$}%
}}}}
\put(7651,-4111){\makebox(0,0)[b]{\smash{{\SetFigFont{8}{9.6}{\familydefault}{\mddefault}{\updefault}{\color[rgb]{0,0,0}$\frac{1}{2h_y}\eey$}%
}}}}
\put(7651,-1711){\makebox(0,0)[b]{\smash{{\SetFigFont{8}{9.6}{\familydefault}{\mddefault}{\updefault}{\color[rgb]{0,0,0}$-\frac{1}{2h_y}\eey$}%
}}}}
\put(6751,-2911){\makebox(0,0)[b]{\smash{{\SetFigFont{8}{9.6}{\familydefault}{\mddefault}{\updefault}{\color[rgb]{0,0,0}$\frac{1}{2h_x}\eex$}%
}}}}
\put(8551,-2911){\makebox(0,0)[b]{\smash{{\SetFigFont{8}{9.6}{\familydefault}{\mddefault}{\updefault}{\color[rgb]{0,0,0}$-\frac{1}{2h_x}\eex$}%
}}}}
\put(7651,-2911){\makebox(0,0)[b]{\smash{{\SetFigFont{8}{9.6}{\familydefault}{\mddefault}{\updefault}{\color[rgb]{0,0,0}${\bf0}$}%
}}}}
\end{picture}%

%% file: gradcartmodif.pstex_t
\begin{picture}(0,0)%
\includegraphics{gradcartmodif.pstex}%
\end{picture}%
\setlength{\unitlength}{3947sp}%
\begingroup\makeatletter\ifx\SetFigFont\undefined%
\gdef\SetFigFont#1#2#3#4#5{%
  \reset@font\fontsize{#1}{#2pt}%
  \fontfamily{#3}\fontseries{#4}\fontshape{#5}%
  \selectfont}%
\fi\endgroup%
\begin{picture}(3120,4014)(6043,-4963)
\put(6751,-4936){\makebox(0,0)[b]{\smash{{\SetFigFont{6}{7.2}{\familydefault}{\mddefault}{\updefault}{\color[rgb]{0,0,0}$i-1$}%
}}}}
\put(7651,-4936){\makebox(0,0)[b]{\smash{{\SetFigFont{6}{7.2}{\familydefault}{\mddefault}{\updefault}{\color[rgb]{0,0,0}$i$}%
}}}}
\put(8551,-4936){\makebox(0,0)[b]{\smash{{\SetFigFont{6}{7.2}{\familydefault}{\mddefault}{\updefault}{\color[rgb]{0,0,0}$i+1$}%
}}}}
\put(6151,-1711){\makebox(0,0)[b]{\smash{{\SetFigFont{6}{7.2}{\familydefault}{\mddefault}{\updefault}{\color[rgb]{0,0,0}$j+1$}%
}}}}
\put(6151,-4111){\makebox(0,0)[b]{\smash{{\SetFigFont{6}{7.2}{\familydefault}{\mddefault}{\updefault}{\color[rgb]{0,0,0}$j-1$}%
}}}}
\put(6151,-2911){\makebox(0,0)[b]{\smash{{\SetFigFont{6}{7.2}{\familydefault}{\mddefault}{\updefault}{\color[rgb]{0,0,0}$j$}%
}}}}
\put(6751,-1711){\makebox(0,0)[b]{\smash{{\SetFigFont{8}{9.6}{\familydefault}{\mddefault}{\updefault}{\color[rgb]{0,0,0}${\bf0}$}%
}}}}
\put(8551,-1711){\makebox(0,0)[b]{\smash{{\SetFigFont{8}{9.6}{\familydefault}{\mddefault}{\updefault}{\color[rgb]{0,0,0}${\bf0}$}%
}}}}
\put(8551,-4111){\makebox(0,0)[b]{\smash{{\SetFigFont{8}{9.6}{\familydefault}{\mddefault}{\updefault}{\color[rgb]{0,0,0}${\bf0}$}%
}}}}
\put(6751,-4111){\makebox(0,0)[b]{\smash{{\SetFigFont{8}{9.6}{\familydefault}{\mddefault}{\updefault}{\color[rgb]{0,0,0}${\bf0}$}%
}}}}
\put(7651,-4561){\makebox(0,0)[b]{\smash{{\SetFigFont{8}{9.6}{\familydefault}{\mddefault}{\updefault}{\color[rgb]{0,0,0}$\frac{1-\sqrt{2}}{2h_y}\eey$}%
}}}}
\put(7651,-2536){\makebox(0,0)[b]{\smash{{\SetFigFont{8}{9.6}{\familydefault}{\mddefault}{\updefault}{\color[rgb]{0,0,0}$-\frac{\sqrt{2}}{h_y}\eey$}%
}}}}
\put(7651,-3361){\makebox(0,0)[b]{\smash{{\SetFigFont{8}{9.6}{\familydefault}{\mddefault}{\updefault}{\color[rgb]{0,0,0}$\frac{\sqrt{2}}{h_y}\eey$}%
}}}}
\put(7426,-2911){\rotatebox{90.0}{\makebox(0,0)[b]{\smash{{\SetFigFont{8}{9.6}{\familydefault}{\mddefault}{\updefault}{\color[rgb]{0,0,0}$\frac{\sqrt{2}}{h_x}\eex$}%
}}}}}
\put(7876,-2911){\rotatebox{270.0}{\makebox(0,0)[b]{\smash{{\SetFigFont{8}{9.6}{\familydefault}{\mddefault}{\updefault}{\color[rgb]{0,0,0}$-\frac{\sqrt{2}}{h_x}\eex$}%
}}}}}
\put(7651,-3736){\makebox(0,0)[b]{\smash{{\SetFigFont{8}{9.6}{\familydefault}{\mddefault}{\updefault}{\color[rgb]{0,0,0}$\frac{1+\sqrt{2}}{2h_y}\eey$}%
}}}}
\put(7426,-4111){\rotatebox{90.0}{\makebox(0,0)[b]{\smash{{\SetFigFont{8}{9.6}{\familydefault}{\mddefault}{\updefault}{\color[rgb]{0,0,0}$\frac{1}{2h_y}\eey$}%
}}}}}
\put(7426,-1711){\rotatebox{90.0}{\makebox(0,0)[b]{\smash{{\SetFigFont{8}{9.6}{\familydefault}{\mddefault}{\updefault}{\color[rgb]{0,0,0}$-\frac{1}{2h_y}\eey$}%
}}}}}
\put(7876,-4111){\rotatebox{270.0}{\makebox(0,0)[b]{\smash{{\SetFigFont{8}{9.6}{\familydefault}{\mddefault}{\updefault}{\color[rgb]{0,0,0}$\frac{1}{2h_y}\eey$}%
}}}}}
\put(7876,-1711){\rotatebox{270.0}{\makebox(0,0)[b]{\smash{{\SetFigFont{8}{9.6}{\familydefault}{\mddefault}{\updefault}{\color[rgb]{0,0,0}$-\frac{1}{2h_y}\eey$}%
}}}}}
\put(8551,-2536){\makebox(0,0)[b]{\smash{{\SetFigFont{8}{9.6}{\familydefault}{\mddefault}{\updefault}{\color[rgb]{0,0,0}$-\frac{1}{2h_x}\eex$}%
}}}}
\put(8551,-3361){\makebox(0,0)[b]{\smash{{\SetFigFont{8}{9.6}{\familydefault}{\mddefault}{\updefault}{\color[rgb]{0,0,0}$-\frac{1}{2h_x}\eex$}%
}}}}
\put(6751,-3361){\makebox(0,0)[b]{\smash{{\SetFigFont{8}{9.6}{\familydefault}{\mddefault}{\updefault}{\color[rgb]{0,0,0}$\frac{1}{2h_x}\eex$}%
}}}}
\put(6751,-2536){\makebox(0,0)[b]{\smash{{\SetFigFont{8}{9.6}{\familydefault}{\mddefault}{\updefault}{\color[rgb]{0,0,0}$\frac{1}{2h_x}\eex$}%
}}}}
\put(8851,-2911){\rotatebox{270.0}{\makebox(0,0)[b]{\smash{{\SetFigFont{8}{9.6}{\familydefault}{\mddefault}{\updefault}{\color[rgb]{0,0,0}$\frac{-1+\sqrt{2}}{2h_x}\eex$}%
}}}}}
\put(6451,-2911){\rotatebox{90.0}{\makebox(0,0)[b]{\smash{{\SetFigFont{8}{9.6}{\familydefault}{\mddefault}{\updefault}{\color[rgb]{0,0,0}$\frac{1-\sqrt{2}}{2h_x}\eex$}%
}}}}}
\put(7651,-1336){\makebox(0,0)[b]{\smash{{\SetFigFont{8}{9.6}{\familydefault}{\mddefault}{\updefault}{\color[rgb]{0,0,0}$\frac{-1+\sqrt{2}}{2h_y}\eey$}%
}}}}
\put(7651,-2161){\makebox(0,0)[b]{\smash{{\SetFigFont{8}{9.6}{\familydefault}{\mddefault}{\updefault}{\color[rgb]{0,0,0}$\frac{-1-\sqrt{2}}{2h_y}\eey$}%
}}}}
\put(6976,-2911){\rotatebox{270.0}{\makebox(0,0)[b]{\smash{{\SetFigFont{8}{9.6}{\familydefault}{\mddefault}{\updefault}{\color[rgb]{0,0,0}$\frac{1+\sqrt{2}}{2h_x}\eex$}%
}}}}}
\put(8326,-2911){\rotatebox{90.0}{\makebox(0,0)[b]{\smash{{\SetFigFont{8}{9.6}{\familydefault}{\mddefault}{\updefault}{\color[rgb]{0,0,0}$\frac{-1-\sqrt{2}}{2h_x}\eex$}%
}}}}}
\end{picture}%